\documentclass[11pt]{amsart}

\usepackage{amssymb,amsmath,bm,amsthm}
\usepackage{graphicx}
\usepackage{multirow}
\usepackage{graphicx}
\usepackage[top=26mm, bottom=26mm, left=23mm, right=25mm]{geometry}
\usepackage[dvips]{color}
\usepackage{bm}

\usepackage{subfigure}
\usepackage{float,enumerate,wrapfig}
\usepackage{algorithm}
\usepackage{algpseudocode}

\newtheorem{remark}{Remark}
\newtheorem{definition}{Definition}
\numberwithin{definition}{section}

\numberwithin{pro}{section}
\newcommand{\mtopp}[1]{{\,\mbox{\footnotesize #1}}}

\newcommand{\Gc}{{\color{black} \Gamma_c}}

\definecolor{Red}{rgb}{1,0,0}
\definecolor{Blue}{rgb}{0,0,1}
\definecolor{Green}{rgb}{0,1,0}
\definecolor{magenta}{rgb}{1,0,.6}
\definecolor{lightblue}{rgb}{0,.5,1}
\definecolor{lightpurple}{rgb}{.6,.4,1}
\definecolor{gold}{rgb}{.6,.5,0}
\definecolor{orange}{rgb}{1,0.4,0}
\definecolor{hotpink}{rgb}{1,0,0.5}
\definecolor{newcolor2}{rgb}{.5,.3,.5}
\definecolor{newcolor}{rgb}{0,.3,1}
\definecolor{newcolor3}{rgb}{1,0,.35}
\definecolor{darkgreen1}{rgb}{0, .35, 0}
\definecolor{darkgreen}{rgb}{0, .6, 0}
\definecolor{darkred}{rgb}{.75,0,0}

\newcommand{\vct}[1]{\bm{\mathsf{#1}}}

\newcommand{\mtx}[1]{\bm{\mathsf{#1}}}

\newcommand{\mtwo}[4]{\left[\begin{array}{cc} #1 & #2 \\ #3 & #4 \end{array}\right]}

\newcommand{\vtwo}[2]{\left[\begin{array}{cc} #1 \\ #2 \end{array}\right]}

\newcommand{\pgnotate}[1]{}

\newcommand{\lsp}{\vspace{3mm}}

\newcommand{\lsm}{\vspace{1mm}}

\newcommand{\bi}{\begin{itemize}}

\newcommand{\ei}{\end{itemize}}

\newcommand{\ben}{\begin{enumerate}}

\newcommand{\een}{\end{enumerate}}

\newcommand{\be}{\begin{equation}}

\newcommand{\ee}{\end{equation}}

\newcommand{\bea}{\begin{eqnarray}}

\newcommand{\eea}{\end{eqnarray}}

\newcommand{\ba}{\begin{align}}

\newcommand{\ea}{\end{align}}

\newcommand{\bse}{\begin{subequations}}

\newcommand{\ese}{\end{subequations}}

\newcommand{\bc}{\begin{center}}

\newcommand{\ec}{\end{center}}

\newcommand{\bfi}{\begin{figure}}

\newcommand{\efi}{\end{figure}}

\newcommand{\bmp}[1]{\begin{minipage}{#1}}

\newcommand{\emp}{\end{minipage}}

\newcommand{\bp}{\begin{proof}}

\newcommand{\ep}{\end{proof}}

\newcommand{\mbf}[1]{{\mathbf #1}}

\begin{document}

\begin{center}
\textbf{A fast direct solver for boundary value problems on 
locally perturbed geometries}

\lsp

\textit{\small Y. Zhang, and A. Gillman\\
Department of Computational and Applied Mathematics, Rice University\\
Yabin.Zhang@rice.edu, adrianna.gillman@rice.edu}

\lsp

\begin{minipage}{0.9\textwidth}\small
\noindent\textbf{Abstract:}
Many applications involve solving several boundary value problems on geometries that 
are local perturbations of an original geometry.  The boundary integral equation for 
a problem on a locally perturbed geometry can be expressed as a low rank update to the 
original system.  A fast direct solver 
 for the new linear system is presented in this paper.  
The solution technique 
utilizes a precomputed fast direct solver for the original geometry to efficiently 
create the low rank factorization of the update matrix and to accelerate the 
application of the Sherman-Morrison formula.  
The method is ideally suited for problems where the local perturbation is the 
same but its placement on the boundary changes and problems where the local 
perturbation is a refined discretization on the same geometry. 
Numerical results illustrate that for fixed local perturbation the method 
is three times faster than building a new fast direct solver from scratch.
\end{minipage}

\end{center}

\section{Introduction}
This manuscript presents a fast direct solver for 
boundary integral equations where the geometry for each problem corresponds 
to a local perturbation of the original geometry.  In particular,
we are interested in problems where the local perturbation to the 
geometry is much smaller than the original geometry.  Since a 
direct solver is constructed, the technique is ideal for problems
with many right hand sides and/or suffer from ill-conditioning due
to geometric complexity.  Boundary value problems involving locally perturbed geometries 
arise in a variety of applications such as optimal design \cite{1989_shape_design}, and 
adaptive discretization techniques \cite{2006_Djokic}. 
For example, finding the optimal placement of an attachment 
to a large geometry which minimizes the radar cross section involves 
solving many problems where the local perturbation is the same but
the placement on the boundary changes.

For many boundary integral equations, the linear
system that results from the discretization of an integral 
equation is amenable to fast direct solvers such as those built from 
 hierarchically semiseparable (HSS) \cite{2007_shiv_sheng,2004_gu_divide}, $\mathcal{H}$-matrix \cite{hackbusch},
hierarchically block separable (HBS) \cite{2012_martinsson_FDS_survey}, hierarchical interpolatory 
factorization (HIF) \cite{2016_ho} and hierarchical off-diagonal low rank (HODLR) \cite{HOLDR} representations.  
These direct solvers utilize the fact that the off-diagonal 
blocks of the dense matrix are low rank. The different
variants correspond to different ways of exploiting this
property.  Let $\mtx{A}$ denote the matrix resulting from discretization
of the boundary integral equation.  The factored approximation of the 
matrix $\mtx{A}$, denoted by $\tilde{\mtx{A}}$, is constructed so that 
$\|\mtx{A}-\tilde{\mtx{A}}\|\leq \epsilon$ for a user defined tolerance $\epsilon$.  
$\tilde{\mtx{A}}$ is called the 
\textit{compressed representation} of $\mtx{A}$.  The 
inverse of the compressed matrix is then constructed via
a variant of a Woodbury formula or by expanding the matrix
out to a larger sparse system and using a sparse direct solver.
We refer the reader to the references for further
details. 

Building from the approach in \cite{2009_martinsson_ACTA},
the solution technique presented in this paper casts the linear 
system for problem on the perturbed geometry as an extended
linear system which consists of a two-by-two block diagonal matrix
plus a low rank update.  The block diagonal matrix has a block 
equal to the matrix for the original geometry.  
By using the Sherman-Morrison formula, the approximate
inverse of the original system can be exploited and the 
approximate inverse of the extended system can be applied rapidly. 
The compressed representation matrix in the original system
is utilized to reduce the cost of computing the low rank factorization
of the update matrix.  The techniques used in this work draw from earlier work in 
\cite{mdirect,BCR,1996_mich_elongated,1994_starr_rokhlin}.

The method presented in this paper is ideally suited for problems where the 
local perturbation is the same over many placements on the boundary or the 
number of removed points on the boundary is not large.  Since the solution 
technique does not modify the original compressed representation, it
can be combined with any fast direct solver.

\subsection{Related work}
 The paper \cite{2016_update} presents a technique for updating 
 the \textit{Hierarchical interpolative factorization} (HIF) 
 of the matrix $\mtx{A}$.  This task involves locating and updating the 
 relevant low rank factors and (potentially) modifying the underlying hierarchical
 tree structure. An approximate inverse is then constructed for
 the updated compressed representation.  The inversion step is one of the most 
 expensive steps in the precomputation of a fast direct solver. 

\subsection{Outline of paper}
This manuscript begins by presenting the boundary integral formulation,
the discretized linear system and the extended linear system for a model 
problem in section \ref{sec:model}.  
The fast direct solver presented in this paper for a locally perturbed
geometry utilizes the factors computed in the fast direct solver for the 
original geometry.  While the method can be utilized in conjunction
with any fast direct solver, we chose to review the HBS method and its 
physical interpretation in section \ref{sec:FDS} for simplicity of presentation.  
Next, the construction of the new fast direct solver is presented in section
\ref{sec:new_fds}.  Finally numerical experiments report on the performance
of the solver in section \ref{sec:numerics}.  Section \ref{sec:summary} 
reviews the methods and highlights the potential of the solution technique.

\section{Model problem}
\label{sec:model}
This section begins by reviewing the boundary integral approach
for solving a Laplace boundary value problem.
Then the technique for writing the linear system corresponding
to the discretized boundary value problem on a geometry that is a local perturbation of the original
is presented in section \ref{sec:local}.

\subsection{Boundary integral equations}
\label{sec:BIE}
\begin{figure}[H]
\begin{center}
\begin{tabular}{c}
 \setlength{\unitlength}{1mm}
\begin{picture}(60,30)(0,0)
\put(15,0){\includegraphics[height=30mm]{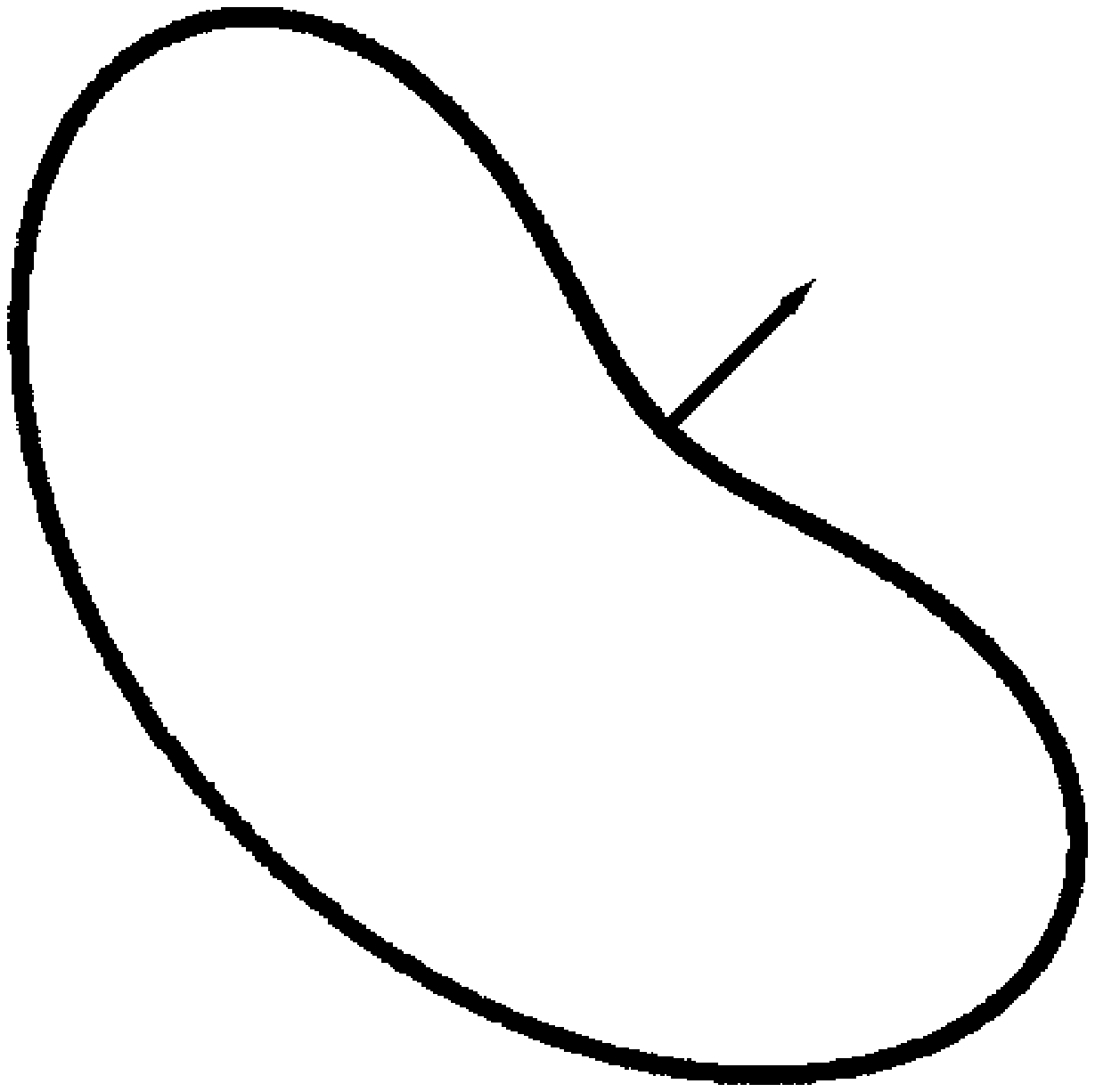}}
 \put(30,10){$\Omega$}
\put(28,28){$\Gamma$}
\put(32,16){$\vct{x}$}
\put(37,19){$\nu_{\vct{x}}$}
\end{picture}\\
(a)
\end{tabular}
\begin{tabular}{c}
 \setlength{\unitlength}{1mm}
\begin{picture}(60,30)(0,0)
\put(15,0){\includegraphics[height=30mm]{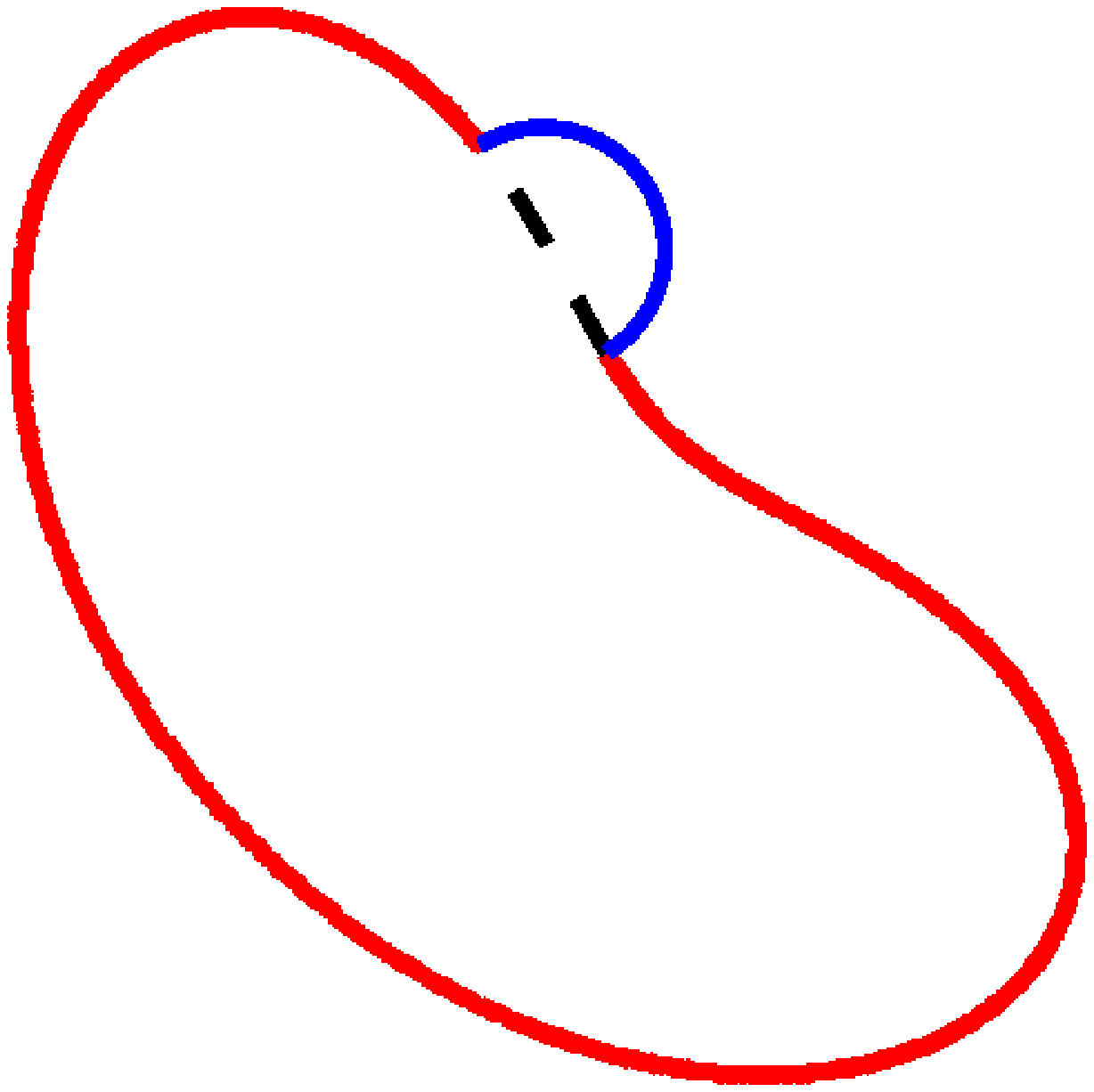}}
 \put(10,20){$\Gamma_{k}$}
\put(30,28){$\Gamma_{p}$}
\put(25,20){$\Gamma_{c}$}
\end{picture}\\
(b)
\end{tabular}
\end{center}
\caption{\label{fig:sample} (a) A sample geometry $\Omega$ with boundary $\Gamma$ and 
outward facing normal vector $\vct{\nu_x}$ at the point $\vct{x}$.  (b) A sample locally perturbed geometry where 
the original boundary is $\Gamma_o = \Gamma_k\cup \Gamma_c$, the portion of the boundary
being removed is $\Gamma_c$, the portion of the original boundary remaining is $\Gamma_k$ 
and the newly added boundary is $\Gamma_p$.}
\end{figure}

For simplicity of presentation, we consider the Laplace boundary 
value problem 
\begin{equation}
\begin{split}
-\Delta u(\vct{x}) =&\ 0, \ \ \quad \vct{x}\in \Omega,\\
u(\vct{x}) =&\ g(\vct{x}),\quad \vct{x} \in  \Gamma.
\end{split}
\label{eq:model}
\end{equation}
Figure \ref{fig:sample}(a) illustrates a sample geometry.  The vector 
$\vct{\nu}_{\vct{y}}$ denotes the outward normal vector at the point $\vct{y}\in\Gamma$.
For $\vct{x}\in \Omega$, we represent the solution to (\ref{eq:model}) 
as a double layer potential
\begin{equation}
 \label{eq:soln}
u(\vct{x}) = \int_{\Gamma}  D(\vct{x},\vct{y})\sigma(\vct{y}) ds(\vct{y}), \qquad \vct{x}\in \Omega,
\end{equation}
where $D(\vct{x},\vct{y}) = \partial_{\vct{\nu}_{\vct{y}}}G(\vct{x},\vct{y})$ 
is double layer kernel, $G(\vct{x},\vct{y}) = -\frac{1}{2\pi}\log|\vct{x}-\vct{y}|$ 
is the fundamental solution and $\sigma(\vct{x})$ is the unknown boundary charge 
distribution.   By taking the limit of $u(\vct{x})$
as $\vct{x}$ goes to the boundary and setting it equal to 
$g(\vct{x})$, we find the boundary charge distribution $\sigma(\vct{x})$
satisfies the following boundary integral equation

\begin{equation}
-\frac{1}{2}\sigma(\vct{x}) +\int_{\Gamma} D(\vct{x},\vct{y})\sigma(\vct{y}) ds(\vct{y}) = g(\vct{x}).
\label{eq:BIE} 
\end{equation}

Discretization of the boundary integral equation (\ref{eq:BIE}) 
with either a Nystr\"om or a boundary element
method results in a linear system of the form 
\begin{equation}\mtx{A}\vct{\sigma} = \vct{g}\label{eq:DBIE}\end{equation}
where the solution $\vct{\sigma}$ is the approximation of 
$\sigma(\vct{x})$ at the discretization points on $\Gamma$.


\subsection{An extended linear system for local perturbations}
\label{sec:local}
In this section, we consider the Laplace boundary 
value problem (\ref{eq:model}) on the geometry 
illustrated in Figure \ref{fig:sample}(b) where 
$\Gamma = \Gamma_{k}\cup\Gamma_p$, $\Gamma_k$ is the remaining portion of 
the original boundary, $\Gamma_c$ is the removed portion of the original 
boundary and $\Gamma_p$ is the 
newly added portion of the boundary.  We define the boundary of the 
original geometry by $\Gamma_o = \Gamma_k\cup\Gamma_c$.

The discretized problem on $\Gamma$ can be expressed as an extended 
linear system \cite{2009_martinsson_ACTA} by
\begin{equation}
\label{eq:extended_system}
\left(\underbrace{
\mtwo{\mtx{A}_{oo}}{\mtx{0}}{\mtx{0}}{\mtx{A}_{pp}}}_{\tilde{\mtx{A}}}
+
\underbrace{\begin{bmatrix}
0 & \begin{pmatrix}-\mtx{A}_{kc}\\-\mtx{B}_{cc}\end{pmatrix} & \mtx{A}_{op}\\
\mtx{A}_{pk} & 0 & 0
\end{bmatrix}}_{\mtx{Q}}
\right)
\begin{pmatrix}
\vct{\sigma}_k\\
\mtx{\sigma}_c\\
\mtx{\sigma}_p
\end{pmatrix}
=\underbrace{\begin{pmatrix}
\vct{f}_k\\ \vct{0}\\ \vct{f}_p
\end{pmatrix}}_{\vct{f}_{\rm ext}}
\end{equation}
where $\mtx{A}_{kc}$ denotes the submatrix of $\mtx{A}_{oo}$ corresponding to
the interaction between $\Gamma_k$ and $\Gamma_c$, $\mtx{A}_{op}$ 
denotes the discretization of the double layer integral operator
on $\Gamma_p$ evaluated on $\Gamma_o$, $\mtx{A}_{pk}$ denotes 
the discretization of the double layer integral operator
on $\Gamma_k$ evaluated on $\Gamma_p$, and 
$\mtx{B}_{cc}$ denotes the sub-matrix of $\mtx{A}_{oo}$ 
corresponding to the interaction of $\Gamma_c$ with itself
but the diagonal entries are set to zero.
The matrix $\mtx{Q}$ is called the \textit{update matrix}.

For many choices of $\Gamma_p$ and $\Gamma_c$, the update matrix $\mtx{Q}$ is low rank.  
These are the cases where section \ref{sec:new_fds} presents a technique for efficiently 
creating the low rank factorization of $\mtx{Q}$ and how to 
exploit the factorization to make a 
fast direct solver.

\section{A fast direct solver for boundary integral equations}
\label{sec:FDS}

In order to construct the low rank factorization 
of $\mtx{Q}$ as efficiently as possible, we reuse factors in the 
compressed representation of $\mtx{A}_{oo}$.  
To give the new work context and introduce vocabulary, this section 
presents a brief review of the construction of HBS representation of 
the matrix $\mtx{A}$ in equation (\ref{eq:DBIE}).  More details 
and the inversion technique are presented in 
\cite{2012_martinsson_FDS_survey} and \cite{2015_bremer_surfaces}.
Other direct solution techniques, such as $\mathcal{H}$-matrix, HSS, etc., 
use similar factorization 
techniques and can be coupled to the new solver.  

Fast direct solvers for the linear system 
in equation (\ref{eq:DBIE}) begin by creating a data-sparse
representation which approximates $\mtx{A}$.  Roughly speaking, 
a data-sparse representation 
of a matrix is a factorization which requires $O(N)$ 
memory to store where $N$ is the number of discretization 
points.  For many fast direct solvers including the HBS 
method, the reduction in memory is achieved by 
exploiting low rank approximations of off-diagonal blocks.
These data-sparse representations also yield fast matrix
vector multiplication and fast inversion schemes.

This section begins by briefly reviewing the Hierarchically block separable
(HBS) representation of a dense matrix in section \ref{sec:HBS}.  Section
\ref{sec:boundary} reviews the physical interpretation of the 
hierarchical method.
Then, section \ref{sec:proxy} presents a fast technique for 
creating the low rank approximations of off-diagonal blocks.

\subsection{HBS representation}
\label{sec:HBS}
This section reviews the HBS representation of a matrix $\mtx{M}$

Consider the $np\times np$ block partition of a matrix 
$\mtx{M}$ into $p\times p$ blocks each of size $n\times n$:
\begin{equation}
\label{eq:yy0}
\mtx{M} \sim \left[\begin{array}{ccccc}
\mtx{D}_{1}   & \mtx{M}_{1,2} & \mtx{M}_{1,3} & \cdots & \mtx{M}_{1,p} \\
\mtx{M}_{2,1} & \mtx{D}_{2}   & \mtx{M}_{2,3} & \cdots & \mtx{M}_{2,p} \\
\vdots    & \vdots    & \vdots    &        & \vdots    \\
\mtx{M}_{p,1} & \mtx{M}_{p,2} & \mtx{M}_{p,3} & \cdots & \mtx{D}_{p}
\end{array}\right]
\end{equation}
Given a desired tolerance $\epsilon$, for each $\tau = 1,\,2,\,\dots,\,p$, 
there exists a constant $k_\tau$ and  $n\times k_\tau$ matrices 
$\mtx{U}_{\tau}$ and $\mtx{V}_{\tau}$ such that each off-diagonal
block $\mtx{M}_{\sigma,\tau}$ of $\mtx{M}$ admits an approximate factorization
$$\|\mtx{U}_{\sigma} \tilde{\mtx{M}}_{\sigma,\tau}\mtx{V}^*_\tau- \mtx{M}_{\sigma,\tau}\|\leq 
\epsilon \qquad \sigma,\tau \in \{1,\,2,\,\dots,\,p\},\quad \sigma \neq \tau.$$
The columns of $\mtx{U}_\sigma$ form a column basis for the columns of all 
off-diagonal blocks in row $\sigma$.  Likewise, the columns of  $\mtx{V}_\tau$ form a 
row basis for all the rows of all off-diagonal blocks in columns $\tau$.  

This factorization allows $\mtx{M}$ to be approximated in the following 
factored form 
\begin{equation}
\label{eq:yy2}
\mtx{M}  \sim \mtx{U}  \tilde{\mtx{M}}  \mtx{V}^{*}  + \mtx{D},\\
\end{equation}
where
$$
\mtx{U} = \mbox{diag}(\mtx{U}_{1},\,\mtx{U}_{2},\,\dots,\,\mtx{U}_{p}),\quad
\mtx{V} = \mbox{diag}(\mtx{V}_{1},\,\mtx{V}_{2},\,\dots,\,\mtx{V}_{p}),\quad
\mtx{D} = \mbox{diag}(\mtx{D}_{1},\,\mtx{D}_{2},\,\dots,\,\mtx{D}_{p}),
$$
and
$$\tilde{\mtx{M}} = \left[\begin{array}{cccc}
0 & \tilde{\mtx{M}}_{12} & \tilde{\mtx{M}}_{13} & \cdots \\
\tilde{\mtx{M}}_{21} & 0 & \tilde{\mtx{M}}_{23} & \cdots \\
\tilde{\mtx{M}}_{31} & \tilde{\mtx{M}}_{32} & 0 & \cdots \\
\vdots & \vdots & \vdots
\end{array}\right].
$$

This is a one level factorization of $\mtx{M}$ which can be 
inverted via a variation of the Sherman-Morrison-Woodbury 
formula (see Lemma 3.1 in \cite{2012_martinsson_FDS_survey} or
\cite{2015_bremer_surfaces}).  

When the matrix $\tilde{\mtx{M}}$ can be factored in the same 
manner, the matrix is called \textit{Hierarchically block separable}
(HBS).  A three level factorization of this kind is expressed as
\begin{equation}
\label{eq:united4}
\mtx{M} \sim \mtx{U}^{(3)}\bigl(\mtx{U}^{(2)}\bigl(\mtx{U}^{(1)}\,\mtx{B}^{(0)}\,
(\mtx{V}^{(1)})^{*} + \mtx{B}^{(1)}\bigr)
(\mtx{V}^{(2)})^{*} + \mtx{B}^{(2)}\bigr)(\mtx{V}^{(3)})^{*} + \mtx{D}^{(3)},
\end{equation}
where the block structure of the factors is 
\begin{equation*}
\begin{array}{cccccccccccccccccc}
\mtx{U}^{(3)} & \mtx{U}^{(2)} & \mtx{U}^{(1)} & \mtx{B}^{(0)} & (\mtx{V}^{(1)})^{*} & \mtx{B}^{(1)} &
(\mtx{V}^{(2)})^{*} & \mtx{B}^{(2)} & (\mtx{V}^{(3)})^{*} & \mtx{D}^{(3)}.\\
\includegraphics[scale=0.3]{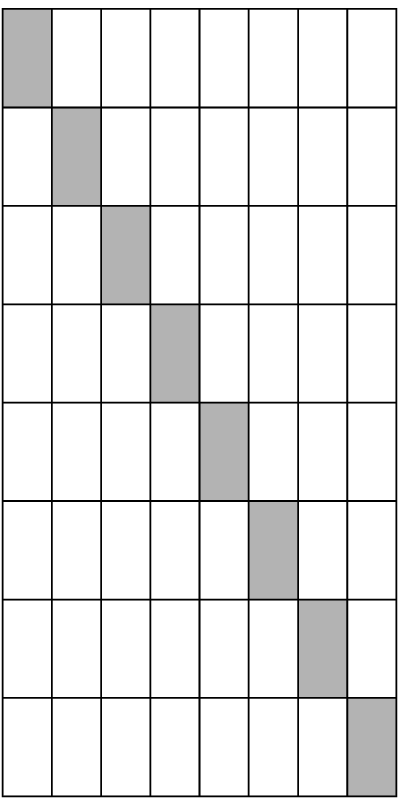}&
\includegraphics[scale=0.3]{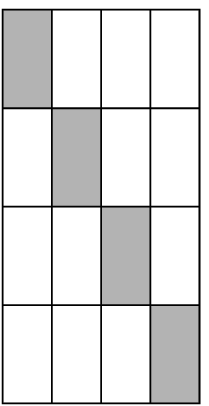}&
\includegraphics[scale=0.3]{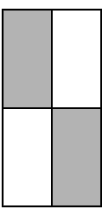}&
\includegraphics[scale=0.3]{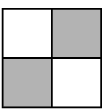}&
\includegraphics[scale=0.3]{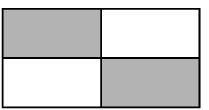}&
\includegraphics[scale=0.3]{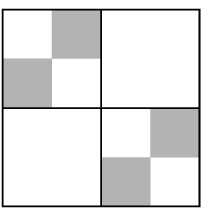}&
\includegraphics[scale=0.3]{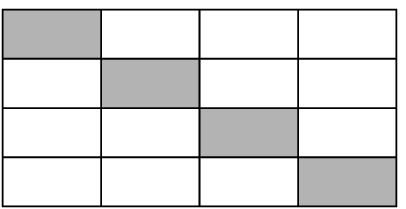}&
\includegraphics[scale=0.3]{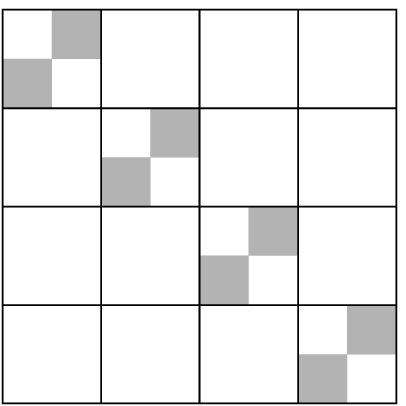}&
\includegraphics[scale=0.3]{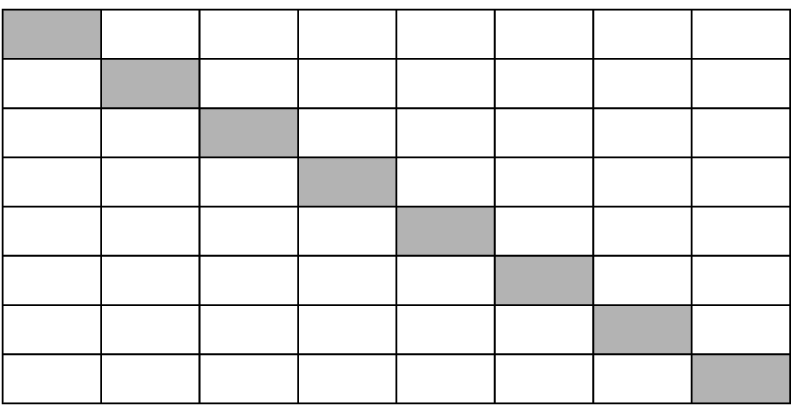}&
\includegraphics[scale=0.3]{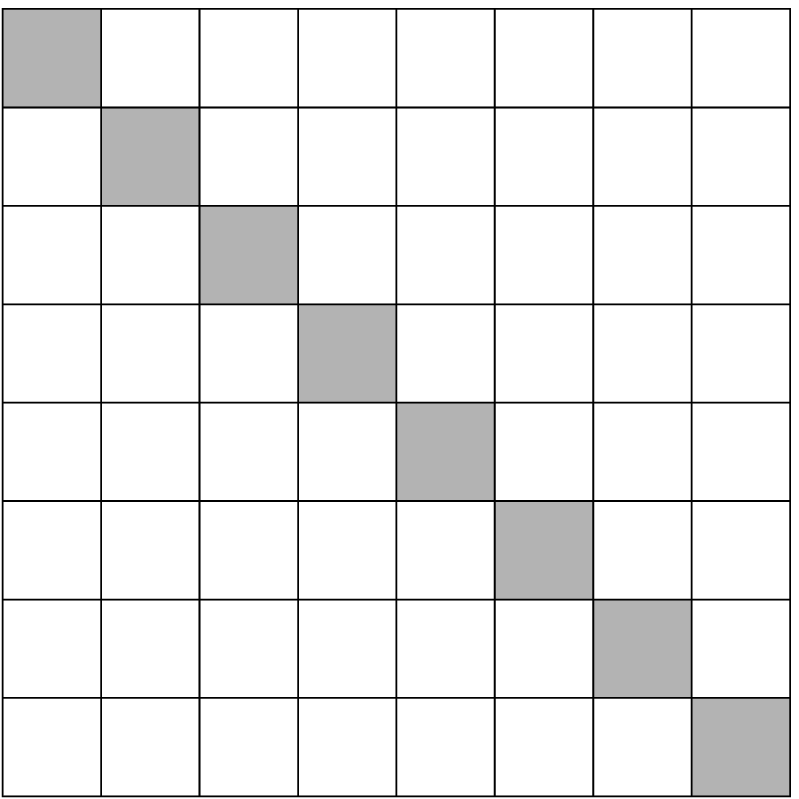}&
\end{array}
\end{equation*}

\subsection{Discretized boundary integral equation to HBS form}
\label{sec:boundary} 

\begin{figure}[H]
\setlength{\unitlength}{1mm}
\begin{picture}(169,41)
\put(20, 0){\includegraphics[height=41mm]{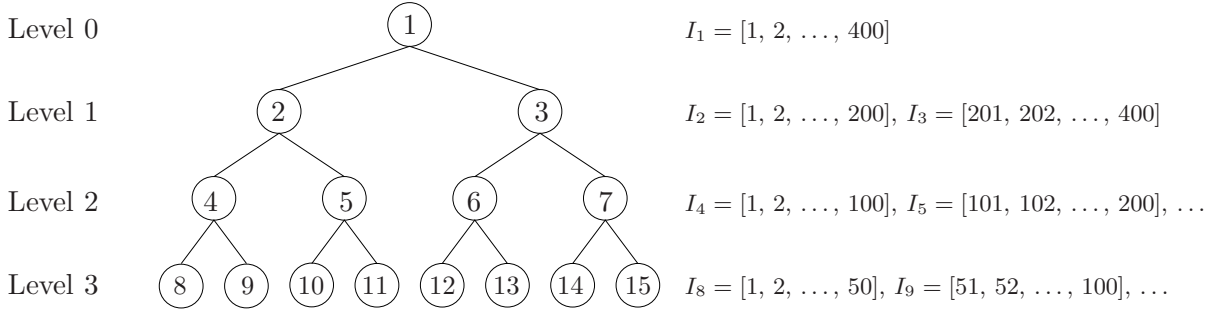}}
\put( 0,36){Level $0$}
\put( 0,25){Level $1$}
\put( 0,13){Level $2$}
\put( 0, 2){Level $3$}
\put(90,36){\footnotesize$I_{1} = [1,\,2,\,\dots,\,400]$}
\put(90,25){\footnotesize$I_{2} = [1,\,2,\,\dots,\,200]$, $I_{3} = [201,\,202,\,\dots,\,400]$}
\put(90,13){\footnotesize$I_{4} = [1,\,2,\,\dots,\,100]$, $I_{5} = [101,\,102,\,\dots,\,200]$, \dots}
\put(90, 2){\footnotesize$I_{8} = [1,\,2,\,\dots,\,50]$, $I_{9} = [51,\,52,\,\dots,\,100]$, \dots}
\put(52.5,36.5){$1$}
\put(35,25){$2$}
\put(70,25){$3$}
\put(26,13){$4$}
\put(44,13){$5$}
\put(61,13){$6$}
\put(78.5,13){$7$}
\put(22, 2){\small$8$}
\put(31, 2){\small$9$}
\put(38.5, 2){\small$10$}
\put(47, 2){\small$11$}
\put(56, 2){\small$12$}
\put(64.5, 2){\small$13$}
\put(73, 2){\small$14$}
\put(82, 2){\small$15$}
\end{picture}
\caption{Numbering of nodes in a fully populated binary tree with $L=3$ levels.
The root is the original index vector $I = I_{1} = [1,\,2,\,\dots,\,400]$.}
\label{fig:tree}
\end{figure}

The HBS representation is based on a binary tree partitioning of an 
index vector $I = [1,\ldots,N]$.  For simplicity, we present the 
technique with a uniform binary tree.  The \textit{root} of the 
tree is $I_1 = I$.  In the next level of the tree, the index vector $I$ is split
into two equilength index vectors $I_2$ and $I_3$. This process is repeated until
each index vector has less than some preset number $n$ of entries.  A \textit{leaf} node
in the tree is an index vector that is not split.  A non-leaf node $\tau$ 
has \textit{children} $\sigma_1$ and $\sigma_2$.  The nodes $\tau$ 
is the \textit{parent} of $\sigma_1$ and $\sigma_2$ if $I_\tau = I_{\sigma_1}\cup I_{\sigma_2}$. Figure \ref{fig:tree} illustrates a three level binary 
tree where $N = 400$ and $n= 50$.

\subsection{Efficient construction of low rank factorizations of off-diagonal blocks}
\label{sec:proxy} 
Constructing the HBS representation of $\mtx{M}$ via pure linear 
algebraic techniques would result in an $O(N^2)$ compression
scheme.  When the matrix $\mtx{M}$ results from the discretization of
a boundary integral equation, physics can be exploited to reduce the 
computational cost of compression to linear.  

First, we note that the index vectors in the binary tree structure
have a physical interpretation.  Since each index corresponds to a 
point on $\Gamma$, an index vector corresponds to a collection of 
discretization points on $\Gamma$ which can be thought of as segment(s) on $\Gamma$.

\begin{figure}[H]
\begin{center}
\begin{tabular}{c}
 \setlength{\unitlength}{1mm}
\begin{picture}(60,50)(0,0)
\put(0,0){\includegraphics[height=50mm]{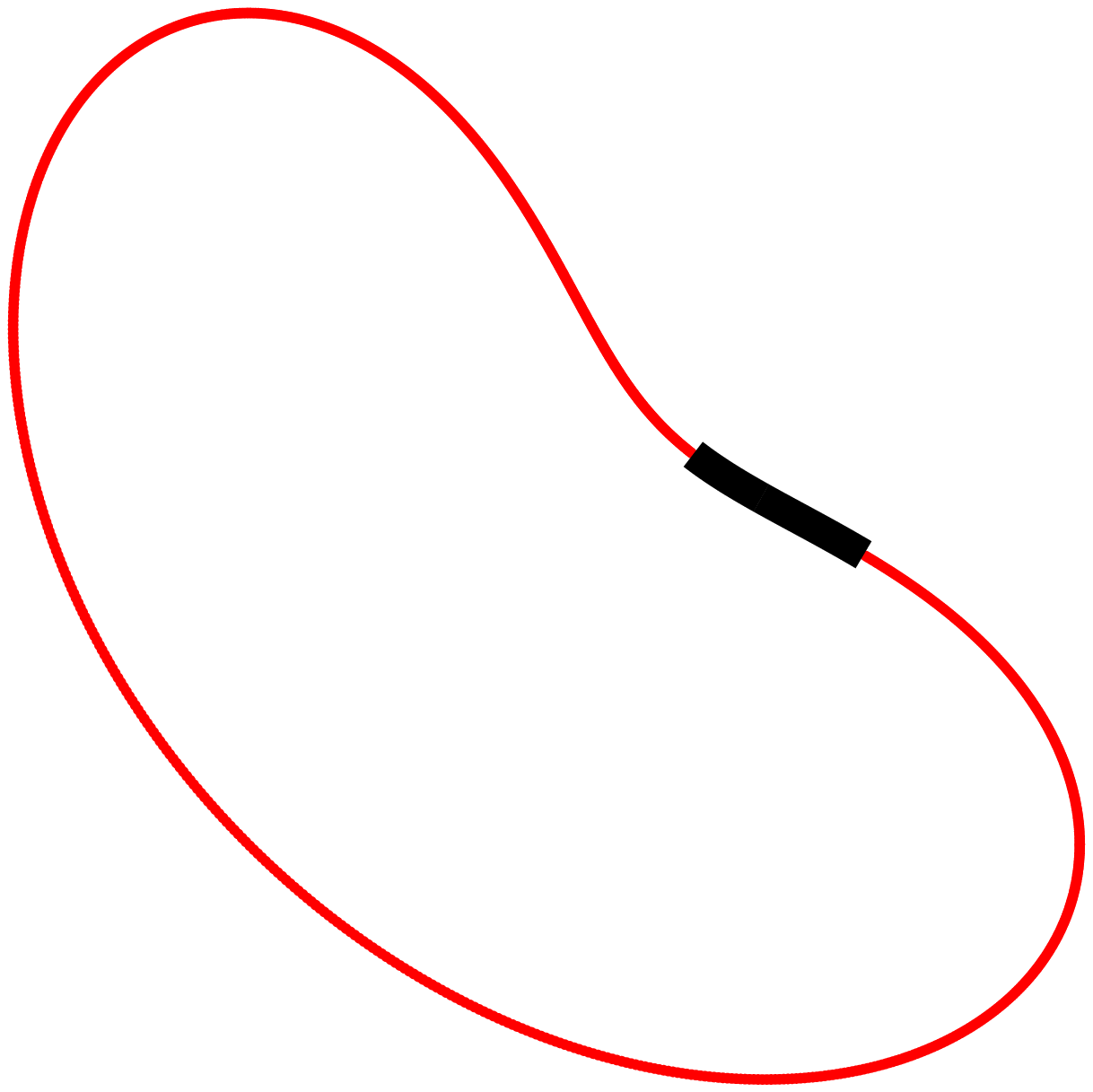}}
 \put(40,23){$\Gamma_{\tau}$}
 \put(08,28){$\Gamma^{c}_\tau$}
\end{picture}\\
(a)
\end{tabular}
\begin{tabular}{c}
 \setlength{\unitlength}{1mm}
\begin{picture}(60,50)(0,0)
\put(0,0){\includegraphics[height=50mm]{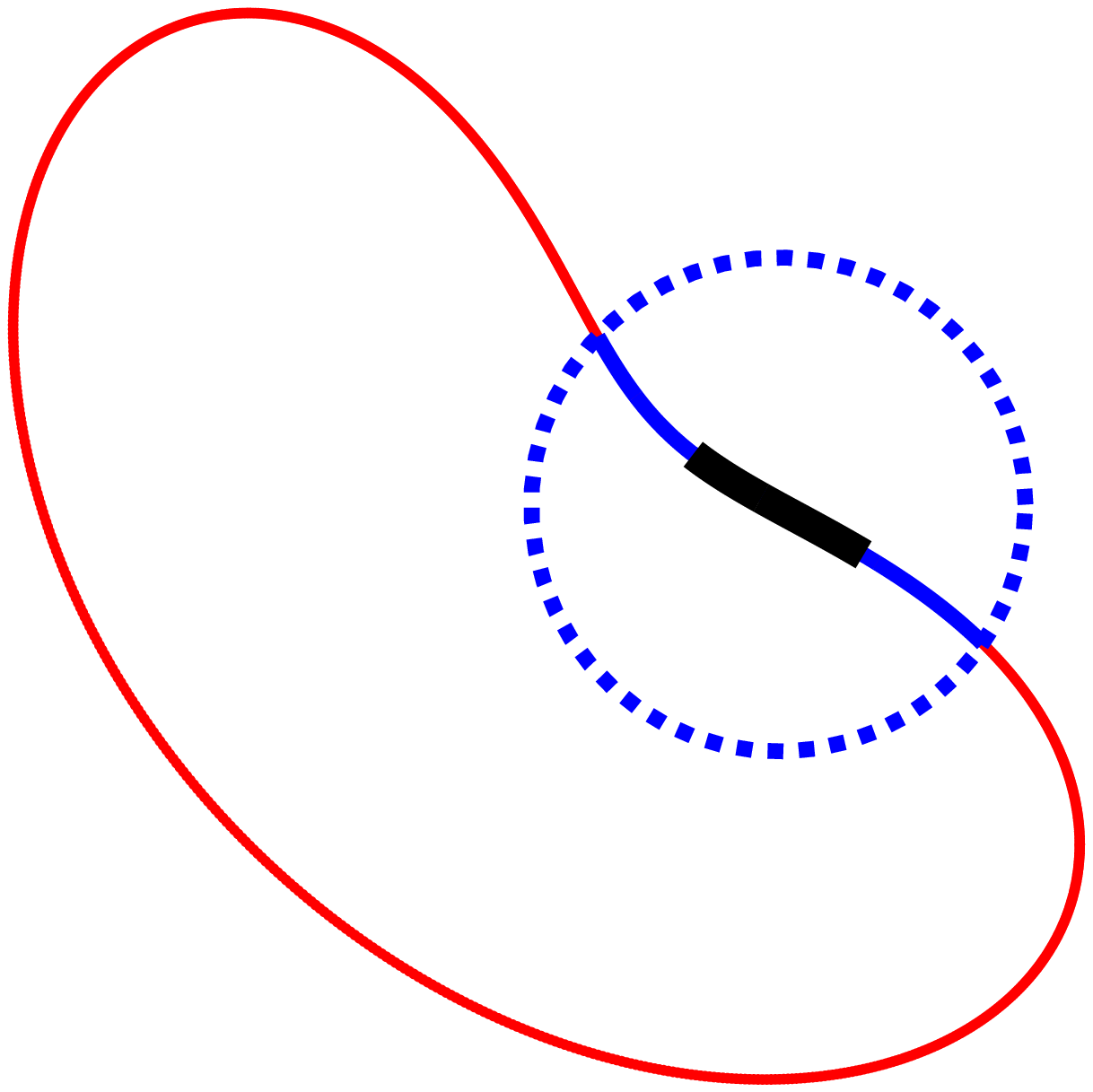}}
 \put(40,23){$\Gamma_{\tau}$}
  \put(52,35){\footnotesize$\color{blue}\Gamma_{\tau}^{\rm{proxy}}$}
 \put(23,07){\footnotesize$\color{red}\Gamma_{\tau}^{\rm{far}}$}
 \put(43,43){\color{blue}\footnotesize$\Gamma_{\tau}^{\rm{near}}$}
 \put(44,42){\color{blue}\vector(-1,-2){5}}
 \put(44,42){\color{blue}\vector(1,-3){5.7}}
\end{picture}\\
(b)
\end{tabular}
\end{center}
\caption{A model geometry with proxy surface. (a) The boundary geometry with $\Gamma_\tau$ in bold line. 
(b) The proxy surface $\Gamma_\tau^{\rm{proxy}}$, drawn in dotted line, that 
separates $\Gamma_\tau^{\rm{far}}$ and $\Gamma_\tau^{\rm{near}}$.}
\label{fig:proxy}
\end{figure}

For presentation purposes, consider the task of creating
the low rank factorization of the sub-matrix of the discretized
linear system (\ref{eq:DBIE}) corresponding to the interaction of $\Gamma_\tau$
and $\Gamma_\tau^c =\Gamma/\Gamma_\tau$ (see Figure \ref{fig:proxy}(a)).  
Let $\mtx{A}_{\tau,c}$ denote this matrix.  

Instead of factoring $\mtx{A}_{\tau,c}$, we partition $\Gamma_\tau^c$ 
into the ``near'' and ``far'' portions.  The \textit{near} portion
of $\Gamma_\tau^c$, denoted by $\Gamma_\tau^{\rm near}$, lies inside a 
\textit{proxy surface} denoted by $\Gamma^{\rm proxy}_{\tau}$.
In this work, we take $\Gamma^{\rm proxy}_\tau$
to be a circle with radius $r\sim 1.5r_{\tau}$, where
$r_\tau$ is the radius of $\Gamma_\tau$, concentric
with $\Gamma_\tau$.  Figure \ref{fig:proxy}(b)
illustrates the proxy surface $\Gamma_\tau^{\rm proxy}$,
proxy points, and near points $\Gamma_\tau^{\rm near}$ 
for $\Gamma_\tau$.  The portion of $\Gamma_\tau^c$ outside of
$\Gamma_\tau^{\rm proxy}$ is the \textit{far} portion of the boundary, denoted by 
$\Gamma_\tau^{\rm far}$.  

From potential theory ideas similar to those 
employed in the fast multipole method \cite{rokhlin1987}, 
it is known that the evaluation of a 
smooth kernel such as (\ref{eq:soln}) for points that are 
far from each other can be expressed with a small number of
basis functions to arbitrary accuracy.  The basis functions
we use are a collection of single poles lying on $\Gamma^{\rm proxy}_\tau$.
This collection of points placed
on the proxy surface are called \textit{proxy points}.
For the examples considered in this paper, we found
that it is enough to have $75$ proxy points.  Let
$\mtx{A}_{\rm proxy}$ denote a matrix characterizing
the interaction between the discretization points on $\Gamma_\tau$
and the proxy points.  Let $\mtx{A}_{\tau,\rm near}$ denote the 
sub-matrix of $\mtx{A}_{\tau,c}$ corresponding to the interaction 
between $\Gamma_{\tau}$ and the near points.   
Then, we compute the low rank factorization of 
$\hat{\mtx{A}} = \left[\mtx{A}_{\tau,\rm near} | \mtx{A}_{\rm proxy}\right]$
which has dimension $N_\tau \times (N_{\rm near}+N_{\rm proxy})$ 
where $N_{\rm near}$ denotes the number of near points, and
$N_{\rm proxy}$ denotes the number of proxy points.  

Instead of computing a QR or SVD, we use
an interpolatory decomposition defined in definition 
\ref{def:interp}.  The underlying algorithms for computing
such a factorization efficiently include rank-revealing QR
\cite{gu1996} and randomized sampling \cite{2007_martinsson_PNAS, 2011_martinsson_randomsurvey}. 
Since
one of the factors is a sub-matrix of the matrix being
factored, the potential theory associated with the 
boundary value problem extends to the factor. Thus
the proxy surface technique can be applied recursively
to create factorizations corresponding to unions
of intervals.  Let 
 \begin{equation}[\mtx{P},J] = {\rm id}(\hat{\mtx{A}})
\label{eq:id}
 \end{equation}
denote the process of computing the interpolatory decomposition where 
$\mtx{P}$ is the interpolation matrix, and $J$ is the corresponding 
index vector.  Algorithm \ref{alg:compress} provides a pseudocode for computing the low rank factorization
of $\mtx{A}_{\tau,c}$ described in this section. 

\begin{definition}
\label{def:interp}
 The \textit{interpolatory decomposition} of a $m\times n$ matrix $\mtx{W}$ that has rank $l$ is
the factorization 
$$ \mtx{W} = \mtx{P}\mtx{W}(J(1:l),:)$$
where $J$ is a vector of integers $j_i$ such $1\leq j_i\leq m$, and $\mtx{P}$ is a $m\times l$ matrix that contains a $l \times l$ identity matrix.
Namely, $\mtx{P}(J(1:l),:) = \mtx{I}_l$.  
\end{definition}

\begin{figure}[ht]
\begin{center}
\fbox{
\begin{minipage}{.9\textwidth}
\begin{center}
\textsc{Algorithm  \ref{alg:compress}} (Efficient factorization of $\mtx{A}_{\tau,c}$)
\end{center}

\lsm

 \textit{Given the boundary $\Gamma$ 
and the boundary segment $\Gamma_\tau$, this 
algorithm computes the low rank factorization of $\mtx{A}_{\tau,c}$ without 
 touching all the entries by using potential theory.
}

\lsm

 \begin{tabbing}
 \hspace{5mm} \= \hspace{5mm} \= \hspace{5mm} \= \hspace{60mm} \= \kill
  
Let $\Gamma^c_\tau = \Gamma/\Gamma_\tau$. \\
  \textit{Make the proxy surface.}\\
Construct $\Gamma^{\rm proxy}$.  \\
  \textit{Find the points near $\Gamma_\tau$.}\\
Set $\Gamma_\tau^{\rm near}$ to be the portion of $\Gamma_\tau^c\subset \Gamma_\tau^{\rm proxy}$.\\
 Let $I_{\rm near}$ denote the discretization points on $\Gamma_\tau^{\rm near}$.\\
  \textit{Make matrices to be factored.}\\
   Let $\mtx{A}_{\rm proxy}$ denote the matrix with the interactions between $\Gamma_\tau$ and $\Gamma^{\rm proxy}$.\\
   Set $\mtx{A}_{\tau,\rm near} = \mtx{A}_\tau(:,I_{\rm near})$.\\
   Set $\hat{\mtx{A}} = [\mtx{A}_{\tau,\rm near}, \ \mtx{A}_{\rm proxy}]$.\\
  \textit{Compute the interpolatory decomposition.}\\
   $[\mtx{P},J] = {\rm id}(\hat{\mtx{A}})$.
  
 \end{tabbing}
\end{minipage}}
\end{center}
\end{figure}

\section{The fast direct solver for the extended linear system}
\label{sec:new_fds}
This section presents the fast direct solver for extended linear
system.  Recall the matrix $\mtx{Q}$ has subblocks $\mtx{A}_{kc}$, $\mtx{A}_{op}$
and $\mtx{A}_{pk}$ which correspond to interactions between small portions
of the boundary and the remainder or portions of the remainder of the boundary.
From section \ref{sec:proxy}, we know these matrices are low rank.  
Let $\mtx{Q}\sim\mtx{L}\mtx{R}$,
where $\mtx{L}\in\mathbb{R}^{n\times k}$ and 
$\mtx{R}\in\mathbb{R}^{k\times n}$, be the rank $k$ approximation
of $\mtx{Q}$.  Then the inverse of $\mtx{A}+\mtx{Q}$ can be 
approximated by the following Sherman-Morrison
formula \cite{golub}
\begin{equation}
\left(\mtx{A}+\mtx{L}\mtx{R}\right)^{-1} = \mtx{A}^{-1} + \mtx{A}^{-1}\mtx{L}\left(\mtx{I}+\mtx{R}\mtx{A}^{-1}\mtx{L}\right)^{-1}\mtx{R}\mtx{A}^{-1}
\label{eq:Woodbury}
\end{equation}
Recall $\mtx{A}$ is a block diagonal matrix with block $\mtx{A}_{oo}$ and $\mtx{A}_{pp}$.
Since $\mtx{A}_{oo}^{-1}$ is already approximated by a fast direct solver and 
the size of $\mtx{A}_{pp}$ is small for applications of interest, the application of 
the approximate inverse of $\mtx{A}$ is fast.  

\begin{remark}
If  $\Gamma_p$ remains the same for multiple perturbed geometries, the cost of constructing $\mtx{A}^{-1}_{pp}$
is not included in the precomputation.  
 In applications where $N_p$ is large, an approximation of $\mtx{A}^{-1}_{pp}$ 
 can be constructed via a fast direct solver. 
\end{remark}

For any local perturbation,
the matrices $\mtx{L}$, $\mtx{R}$, $\mtx{A}^{-1}\mtx{L}$ 
and $\left(\mtx{I}+\mtx{R}\mtx{A}^{-1}\mtx{L}\right)^{-1}$ 
need only be computed once.  The construction of the factorizations and 
$\left(\mtx{I}+\mtx{R}\mtx{A}^{-1}\mtx{L}\right)^{-1}$ comprise
the \textit{precomputation} of the new solver.  Once constructed,
the application of the Woodbury formula (\ref{eq:Woodbury})
can be evaluated for linear cost with small constant.  

In order for this
to be a fast direct solver, the low rank factorization of $\mtx{Q}$ must
scale linearly with the number of discretization points on the original 
geometry.  Before detailing how to efficiently factorize $\mtx{Q}$,
 we introduce some notation.  
Let $N_{k}$ denote the number of discretization points on $\Gamma_k$, $N_{c}$ 
denote the number of discretization points on $\Gamma_c$, and $N_p$ denote the
number of discretization points on $\Gamma_p$.  Then the number of discretization 
points on the original boundary $\Gamma_o$ is $N_o = N_c+N_k$, the number
of discretization points on the new geometry is $N_n = N_k+N_p$, and 
the dimension of $\mtx{Q}$ is $N_{\rm ext}\times N_{\rm ext}$ where 
$N_{\rm ext} = N_o+N_p$.  

The remainder of this section describes how to efficiently construct the 
low rank factorization of $\mtx{Q}$.  The factorization is achieved by 
constructing low rank factorizations of the submatrices.  Let 
\begin{equation}
\begin{array}{cccccccccc}
\mbf{A}_{kc} & \approx & \mbf{L}_{kc} &\mbf{R}_{kc}, &      & &\mbf{A}_{op} &  \approx & \mbf{L}_{op} &\mbf{R}_{op}, \mbox{ and }\\
N_k \times N_c & & N_k\times k_{kc} &  k_{kc}\times N_c &    && N_o \times N_p & & N_o\times k_{op} &  k_{op}\times N_p\\
& & &\mbf{A}_{pk} &  \approx & \mbf{L}_{pk} &\mbf{R}_{pk} &&&\\
& & &N_p \times N_k & & N_p\times k_{pk} &  k_{pk}\times N_k&&&\\
\end{array}
\end{equation}
denote the low rank factorizations of the submatrices.  Then 
\begin{equation}
\label{equ:QeqLR}
\begin{array}{cccc}
\mathbf{Q} &= &\mathbf{L} & \mathbf{R}\\
N_{\rm ext}\times N_{\rm ext}& & N_{\rm ext} \times k & k \times N_{\rm ext}
\end{array}
\end{equation}
where
\begin{equation*}
\mathbf{L} =\begin{bmatrix}\begin{pmatrix} -\mathbf{L}_{kc} &\mathbf{0}\\ 0& -\mathbf{B}_{cc} \end{pmatrix} & \mathbf{L}_{op}\\
\mathbf{L}_{pk} & \mathbf{0} \\
\end{bmatrix}  
\mbox{ and }
\mathbf{R}=
\begin{bmatrix}
\mathbf{R}_{pk} & \mathbf{0} &\mathbf{0}\\
\mathbf{0} &\begin{pmatrix} \mathbf{R}_{kc} \\\mathbf{I}_{cc} \end{pmatrix} & \mathbf{0}\\
\mathbf{0} & \mathbf{0} & \mathbf{R}_{op}\\
\end{bmatrix}.
\end{equation*}

Section \ref{sec:LR_kc} and section \ref{sec:LR_op} present a linearly scaling 
technique for factorizing $\mtx{A}_{kc}$ and $\mtx{A}_{op}$ respectively.  
The low rank factorization of $\mtx{A}_{pk}$ is achieved via a technique 
similar to that as factoring $\mtx{A}_{op}$ and thus is not presented.

\subsection{The efficient factorization of $\mtx{A}_{kc}$}
\label{sec:LR_kc}
The matrix $\mtx{A}_{kc}$ is a sub-matrix of $\mtx{A}_{oo}$ and
thus much of the low rank factorization can be extracted from the 
 compressed representation of $\mtx{A}_{oo}$.  This section 
presents a technique for efficiently creating a low rank factorization
of $\mtx{A}_{kc}$ by reusing as much information as possible from 
the compressed representation of $\mtx{A}_{oo}$.  Roughly speaking, many of
the factors can be found by collecting the $\mtx{U}$ or $\mtx{V}$ 
factors in the HBS factorization (such as (\ref{eq:united4})).  For 
simplicity of presentation, this section is limited to collecting the 
$\mtx{U}$ factors.  

The algorithm begins by first constructing the low rank factorizations
that cannot be extracted from the HBS representation.  This 
consists of leaf boxes $\tau$ that have points in $\Gamma_k$ and $\Gamma_c$.  
Let $W$ denote the set of all such leaf boxes.
For each $\tau\in W$, the points on $\tau\in\Gamma_k$ are identified 
and labeled $J_{\tau,k}$ and the matrix  
$\mtx{A}_{kc}(J_{\tau,k},:)$ is compressed via the method presented in 
section \ref{sec:proxy}.  The result is the interpolation matrix $\mtx{P}_\tau$ and 
index vector $J_\tau$.  

Next we extract as much of the information from the HBS representation
of the matrix $\mtx{A}_{oo}$ as possible.  Let $T$ denote the set of 
boxes $\tau\subset \Gamma_k$.  Information is extracted by going
through the boxes in $T$ from smallest to the largest.     
 For a leaf box $\tau\in T$, 
the interpolation matrix $\mtx{U}_\tau$ is extracted from the factorization.  For a non-leaf
box $\tau\in T$, let $\sigma_1$ and $\sigma_2$ denote the children boxes.  Then
let $\mtx{P}_\tau$ denote the interpolation matrix extracted from the HBS 
representation.  Then 
$$\mtx{U}_{\tau} = \mtwo{\mtx{U}_{\sigma_1}}{\mtx{0}}{\mtx{0}}{\mtx{U}_{\sigma_2}}\mtx{P}_\tau.$$

\begin{figure}[H]
\begin{center}
\begin{picture}(100,130)(01,01)
\put(-50,-10){\includegraphics[scale=0.35]{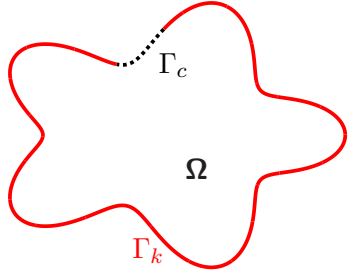}}
\put(60,50){$\mtx{\Omega}$}
\put(40,20){${\color{red}\Gamma_k}$}
\put(50,90){${\color{black}\Gamma_c}$}
\end{picture}
\caption{Geometry used to test the approximate factorization techniques
in section \ref{sec:LR_kc} for $\mtx{A}_{kc}$.  The solid red line 
denotes $\Gamma_k$ and the black dotted line denotes $\Gc$.  The approximate ranks are 
reported in Table \ref{tab:A_kc_rank_demo}. }
\label{fig:A_kc_rank_demo}
\end{center}
\end{figure}

Let $V = W\cup T = \{\tau_1,\ldots,\tau_m\}$.   A low rank factorization can result from letting 
 $\mtx{L}_{kc}$ denote the block diagonal matrix where the subblocks are the matrices 
$\mtx{U}_{\tau_j}$ for $j =1,\ldots,m$ and $\mtx{R}_{kc} = \mtx{A}_{kc}(J,:)$ 
where ${J} = [{J}_{\tau_1},\ldots,{J}_{\tau_m}]$ denotes the corresponding
list of indices.  Unfortunately, 
the size of these factors is significantly larger than optimal.  To prevent this
from hampering the performance of the solver, we must do an additional compression 
step. We call this extra step \textit{recompression}.  Table \ref{tab:A_kc_rank_demo} 
illustrates the approximate ranks computed via the different compression schemes for
the geometry illustrated in Figure \ref{fig:A_kc_rank_demo}, the number of points 
$N_k$ on $\Gamma_k$, the number of points $N_c$ on $\Gamma_c$,
the length $k^0$ of the index vector ${J}$, the size $k$ of the second dimension of 
$\mtx{L}_{kc}$ after recompression, and the optimal rank $k^{\rm opt}$ of $\mtx{A}_{kc}$.  
The size of the factors is close to optimal after the recompression step.

To \textit{recompress} the factorization, we focus our attention on the submatrices
of $\mtx{A}_{kc}$ which contain all the rank information.  To construct the 
low rank factorization of $\mtx{A}_{kc}(J,:)$, we start by constructing
the low rank factorization of $\mtx{A}_{kc}([{J}_{\tau_1},{J}_{\tau_2}],:)$ via the 
method in section \ref{sec:proxy}.  The result is an interpolation matrix $\mtx{P}_{12}$ and 
index vector ${J}_{12}$.  We proceed by constructing the low rank factorization 
of $\mtx{A}_{kc}([{J}_{12},{J}_{\tau_3}],:)$ via the method in section \ref{sec:proxy}.
This continues until all $m$ boxes have been processed.  The new matrix $\mtx{L}_{kc}$ is the 
old matrix multiplying by the interpolation matrices and the final index vector defines $\mtx{R}_{kc}$ as
 $\mtx{R}_{kc}=\mtx{A}_{kc}(J,:)$.  Let $k_{kc}$ denote the length of the index vector $J$. 
 Table \ref{tab:A_kc_rank_demo} illustrates the near optimal rank $k_{kc}$ resulting from this procedure.
Algorithm \ref{alg:recompress} provides a pseudocode for the recompression scheme.

Algorithm \ref{alg:A_kc} provides a pseudocode for the algorithm presented in this section.

\begin{table}[H]
\begin{center}
\begin{tabular}{@{} | c | c | c | c  | c | @{}} 
 \hline
$N_k$  & $N_c$ & $k^0$ & $k$ & $k_\mtopp{opt}$ \\
\hline
1200&  80 & 245 & 17 & 15\\
4800 & 320 & 319 & 17 & 15\\
19200 & 1280 & 361 & 17 & 15\\
\hline
   \end{tabular}
   \caption{The rank numbers of $\mtx{A}_{kc}$ of the test geometry shown 
   in Figure \ref{fig:A_kc_rank_demo}. $k^0$ denotes the length of the
   index vector $J$ produced by the factorization technique in Algorithm 
   \ref{alg:A_kc} prior to the recompression step, and $k$ denotes the 
   length of the index vector $J$  after the recompression step.  $k_\mtopp{opt}$
   is the number of singular values of $\mtx{A}_{kc}$ that are greater 
   than $\epsilon=1\times 10^{-10}$ and is considered as the optimal rank numbers for compression. }
\label{tab:A_kc_rank_demo}
\end{center}
\end{table}

\begin{figure}[ht]
\begin{center}
\fbox{
\begin{minipage}{.9\textwidth}
\begin{center}
\textsc{Algorithm  \ref{alg:A_kc}} (Efficient factorization of $\mtx{A}_{kc}$)
\end{center}

\lsm

\textit{Given the HBS representation of $\mtx{A}_{oo}$, and the portions of the 
boundary $\Gamma_c$ and $\Gamma_k$, this algorithm computes the low rank factorization
of $\mtx{A}_{kc}$ which is a sub-matrix of $\mtx{A}_{oo}$ by reusing as much of the
precomputed HBS representation as possible.
}

\lsm

 \begin{tabbing}
 \hspace{5mm} \= \hspace{5mm} \= \hspace{5mm} \= \hspace{60mm} \= \kill

 \textit{Factor leaf boxes on $\Gamma_k$ that have points on $\Gamma_c$.}\\
 Let $W$ denote the set of leaf boxes factored in this loop. \\
 \textbf{loop} over leaf boxes $\tau$\\
 \> \textbf{if} $\tau\cap \Gamma_c \neq \emptyset$ and $\tau\cap \Gamma_k \neq \emptyset$,\\
\> \> Let ${I}_k$ denote the indices of points in $\tau\cap\Gamma_k$.\\
\> \> Let $\mtx{U}_{\tau_k}$ and ${J}_{\tau_k}$ be the interpolation matrix and 
index vector resulting from \\
\> \> \> applying Algorithm \ref{alg:compress} to $\mtx{A}_{kc}({I}_k,:)$.\\
 \> \textbf{end if}\\
 \textbf{end loop}\\
\\
 
\textit{Extract other factors from HBS representation of $\mtx{A}_{oo}$.}\\
 Let $T$ denote the set containing the largest boxes touched in the loop.\\
\textbf{loop} over boxes $\tau$ from smallest to largest, \\
\> \textbf{if} $\tau\subset \Gamma_k$, \\
\>\> \textbf{if} $\tau$ is a leaf box \\
\> \> \> Set $\mtx{U}_\tau = \mtx{P}_\tau$ where $\mtx{P}_\tau$ is the interpolation matrix 
      from the HBS representation \\
 \>\>\> \hspace{.2cm}     of $\mtx{A}_{oo}$.\\
     \> \> \textbf{else} \\
\> \> \> Let $\sigma_1$ and $\sigma_2$ denote the children of $\tau$.\\
\> \> \> $\mtx{U}_{\tau} = \mtwo{\mtx{U}_{\sigma_1}}{\mtx{0}}{\mtx{0}}{\mtx{U}_{\sigma_2}}\mtx{P}_\tau$.\\
\> \> \textbf{end if}\\
\> \textbf{end if}\\
 \textbf{end loop}\\
 \\
 
 Let $V = W\cup T = \{\tau_1,\ldots,\tau_m\}$.\\
 Let $\mtx{L}_{kc}$ denotes the block diagonal matrix with the $\mtx{U}_\tau$ matrices from $V$ as the blocks.\\

 \\

    Use algorithm \ref{alg:recompress} to recompress the factors. Then\\
     $\mtx{L}_{kc} = \mtx{L}_{kc}\mtx{L}$ 
     and 
   $\mtx{R}_{kc} = \mtx{A}_{kc}(J,:)$.
 \end{tabbing}
\end{minipage}}
\end{center}
\end{figure}

\begin{figure}[ht]
\begin{center}
\fbox{
\begin{minipage}{.9\textwidth}
\begin{center}
\textsc{Algorithm  \ref{alg:recompress}} (Recompression scheme to remove extra degrees of freedom)
\end{center}

\lsm

\textit{Given a set of boxes $V = \{\tau_1,\ldots,\tau_m\}$, the corresponding
indices $J_{\tau_1},\ldots,J_{\tau_m}$ and the original matrix $\mtx{M}$, this algorithm 
efficiently creates a low rank factorization of $\mtx{M}(K,:)$ where 
$K = \{J_{\tau_1},\ldots,J_{\tau_m}\}$.  
}

\lsm

 \begin{tabbing}
  \hspace{5mm} \= \hspace{5mm} \= \hspace{5mm} \= \hspace{60mm} \= \kill

  Let ${J} = [{J}_{\tau_1},{J}_{\tau_2}]$ denote the indices for $\tau_1$ and $\tau_2$.\\
  $[\mtx{P},{J}] = {\rm compress}(\mtx{M}(J,:))$.\\
  $\mtx{L}= \mtwo{\mtx{P}}{0}{0}{\mtx{I}}$.\\ 
  \textbf{loop} over remaining boxes $\tau_j$ in $V$, \\
   \> ${J} = [{J},{J}_{\tau_j}]$.\\
Let $\mtx{P}$ and ${J}$ be the interpolation matrix and 
index vector resulting from \\
\> \> applying Algorithm \ref{alg:compress} to $\mtx{M}(J,:)$.\\
   \> $\mtx{L} = \mtx{L}\mtwo{\mtx{P}}{0}{0}{\mtx{I}}$.\\
   \textbf{end loop}\\
 \end{tabbing}
\end{minipage}}
\end{center}
\end{figure}

\subsection{The efficient compression of $\mtx{A}_{op}$}
\label{sec:LR_op}

This section presents an efficient technique for constructing
the low rank factorization of $\mtx{A}_{op}$.  While this matrix
is not a sub-matrix of $\mtx{A}_{oo}$, the far field information 
from the HBS representation of $\mtx{A}_{oo}$ can be reused in 
constructing the low rank factorization of $\mtx{A}_{op}$. 
Recall from section \ref{sec:proxy} that all far field interactions
can be captured via the interaction with a proxy surface.  This
means we can reuse the interpolation matrices for all boxes that 
are far from $\Gamma_p$.  

The technique for creating this factorization is based on partitioning
$\Gamma_o$ into the portions near $\Gamma^{\rm near}$ and far $\Gamma^{\rm far}$ 
from $\Gamma_p$.  For all boxes $\tau$ contained in $\Gamma^{\rm far}$, the technique for extracting
the factors from the HBS representation of $\mtx{A}_{oo}$ is the same as in
section \ref{sec:LR_kc}.  To compress near field interactions, the technique
is similar to creating an HBS factorization from scratch.

The compression of the near field begins by creating a new binary tree (similar to the 
one in Figure \ref{fig:tree}) for the points $I_{\rm near}$ on $\Gamma^{\rm near}$.  
Let $\mathcal{T}_{near}$ denote the binary tree, then the factorization  
$\mtx{A}(I_{near},I_p)$ can be constructed via a nested factorization.  If a box $\tau$ 
is far from $\Gamma_p$, a proxy surface is used to create the low rank factorization of 
$\mtx{A}_{\tau,p}$.  If a box $\tau$ is near $\Gamma_p$, a proxy surface $\Gamma^{\rm proxy}$ 
is placed around $\tau$ and any points in $\Gamma_p$ that are inside $\Gamma^{\rm proxy}$ are
near $\tau$.  Let $I_p$ denote the indices of those points.  The interpolatory decomposition 
is then applied to $[\mtx{A}_{op}(\tau,I_p), \mtx{A}^{\rm proxy}]$ where $\mtx{A}^{\rm proxy}$ 
denotes the interaction between $\tau$ and the proxy surface.  
Let $\mtx{P}$ and $J$ denote the interpolation matrix and index vector respectively resulting
from applying the interpolatory decomposition to $[\mtx{A}_{op}(\tau,I_p), \mtx{A}^{\rm proxy}]$. 
Algorithm \ref{alg:geo} provides a pseudocode for factorizing a near field interaction.

\begin{figure}[ht]
\begin{center}
\fbox{
\begin{minipage}{.9\textwidth}
\begin{center}
\textsc{Algorithm  \ref{alg:geo}} (Efficient factorization of $\mtx{A}_{\tau,p}$)
\end{center}

\lsm

 \textit{Given $\Gamma_\tau$, and $\Gamma_p$, this
algorithm computes the low rank factorization of $\mtx{A}_{\tau,p}$ without 
 touching all the entries by using potential theory.
}

\lsm

 \begin{tabbing}
 \hspace{5mm} \= \hspace{5mm} \= \hspace{5mm} \= \hspace{60mm} \= \kill
  
  \textit{Make the proxy surface.}\\
Construct $\Gamma^{\rm proxy}$.  \\
  \textit{Find the points on $\Gamma_p$ near $\Gamma_\tau$.}\\
Set $\Gamma^{\rm near}$ to be the portion of $\Gamma_p\subset \Gamma^{\rm proxy}$.\\
 Let $I_{\rm near}$ denote the discretization points on $\Gamma^{\rm near}$.\\
  \textit{Make matrices to be factored.}\\
   Let $\mtx{A}_{\rm proxy}$ denote the matrix with the interactions between $\Gamma_\tau$ and $\Gamma^{\rm proxy}$.\\
   Set $\mtx{A}_{\tau,\rm near} = \mtx{A}_{\tau,p}(:,I_{\rm near})$.\\
   Set $\hat{\mtx{A}} = [\mtx{A}_{\tau,\rm near} \ \mtx{A}_{\rm proxy}]$.\\
  \textit{Compute the interpolatory decomposition.}\\
   $[\mtx{P},J] = {\rm id}(\hat{\mtx{A}})$.
 \end{tabbing}
\end{minipage}}
\end{center}
\end{figure}

Let $\mtx{L}_{\rm far}$ and $\mtx{R}_{\rm far} = \mtx{A}_{op}(J_{\rm far},:)$ denote
the low rank factorization of the far field where $J_{\rm far}$ denotes the 
index vector resulting from the far field compression.  Likewise let $\mtx{L}_{\rm near}$ and
$\mtx{R}_{\rm near} = \mtx{A}_{op}(J_{\rm near},:)$ denote the low rank factorization of the near
field.  A low rank factorization of $\mtx{A}_{op}$ can be $\mtx{L}_{op}\mtx{R}_{op}$ where  
$\mtx{L}_{op} = \mtwo{\mtx{L}_{\rm far}}{0}{0}{\mtx{L}_{\rm near}}$ and 
$\mtx{R}_{op} = \vtwo{\mtx{R}_{\rm far}}{\mtx{R}_{\rm near}}$.
As seen in creating the low rank factorization of $\mtx{A}_{kc}$, the approximate rank given
by this factorization is likely far from optimal.  Thus another compression is necessary.  Let 
$J_{\rm tot} = [J_{\rm far}, J_{\rm near}]$.  Applying the interpolatory decomposition
to $\mtx{A}_{op}(J_{\rm tot},:)$ results in another interpolation matrix $\mtx{P}$ and 
index vector $J$.  Thus the final factors are 
$\mtx{L}_{op} = \mtx{L}_{op}\mtx{P}$ and $\mtx{R}_{op} = \mtx{A}_{op}(J,:)$ where the 
approximate rank $k_{op}$ is the length of the index vector $J$.

Algorithm \ref{alg:A_op} gives a pseudocode for the technique presented in this section.

\begin{figure}[H]
\begin{center}
\begin{tabular}{lcr}
\begin{picture}(100,120)(01,01)
\put(-70,-0){\includegraphics[scale=0.35]{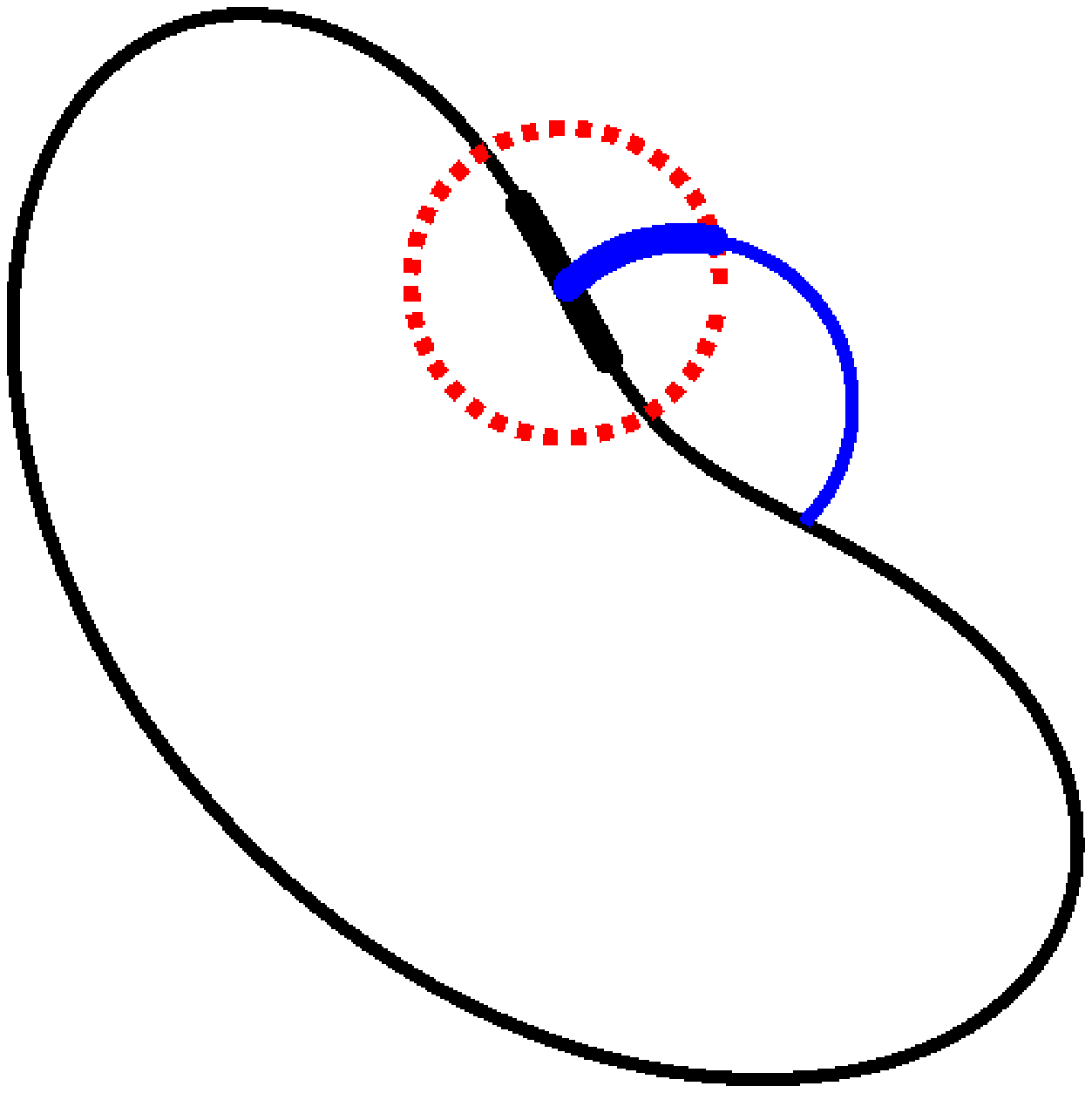}}
\put(00,40){$\mtx{\Omega}$}
\put(-45,20){${\Gamma_k}$}
\put(30,95){${\color{blue}\Gamma_p}$}
\put(-20,70){${\color{red}\Gamma^{\rm proxy}}$}
\put(-10,90){$\Gamma_\tau$}
\end{picture} & {\mbox{\hspace{0.2cm}}} & 

\begin{picture}(100,150)(01,01)
\put(5,-0){\includegraphics[scale=0.35]{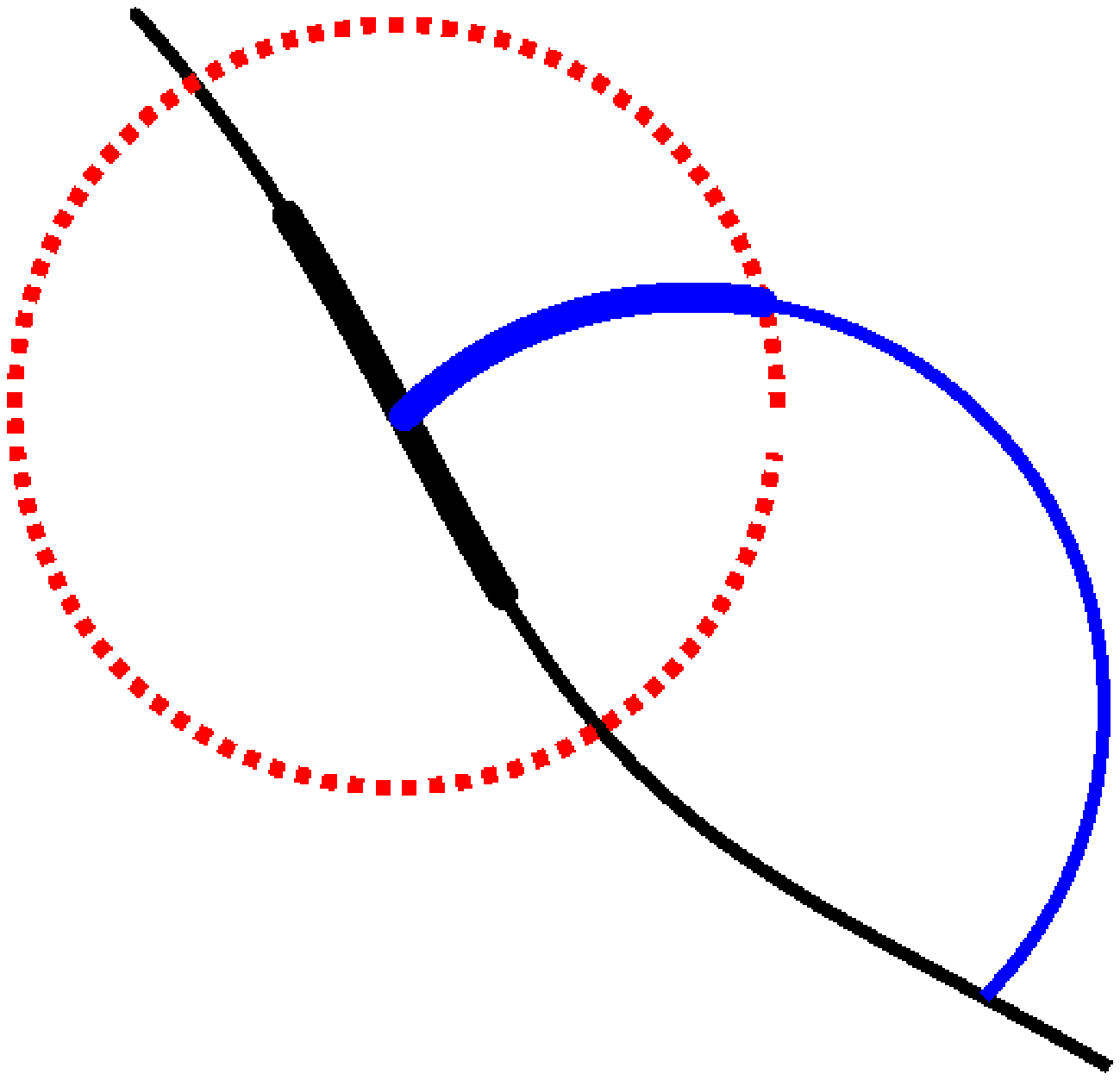}}
\put(50,65){$\Gamma_\tau$}
\put(40,30){${\color{red}\Gamma_{\rm proxy}}$}
\put(65,80){${\color{blue}\Gamma_{\rm near}}$}
\end{picture}\\
(a) & & (b)
\end{tabular}
\caption{Illustration of a geometry for compressing $\mtx{A}_{op}$ when a box $\tau$ 
is near $\Gamma_p$. The whole geometry, $\Gamma_\tau$, and $\Gamma^{\rm proxy}$ are illustrated
in (a).  A close up including $\Gamma^{\rm near}$ is illustrated in (b).}
\label{fig:near_field}
\end{center}
\end{figure}

\begin{figure}[hp]
\begin{center}
\fbox{
\begin{minipage}{.9\textwidth}
\begin{center}
\textsc{Algorithm  \ref{alg:A_op}} (Efficient factorization of $\mtx{A}_{op}$)
\end{center}

\lsm

\textit{Given the HBS representation of $\mtx{A}_{oo}$, $\Gamma_o$ partitioned
into the portion near $\Gamma^{\rm near}$ and far $\Gamma^{\rm far}$ from $\Gamma_p$
and a binary tree $\mathcal{T}_{\rm near}$ defined for points on $\Gamma^{near}$, 
this algorithm computes the low rank factorization
of $\mtx{A}_{op}$ reusing as much of the HBS representation of $\mtx{A}_{oo}$ as 
possible. 
}

\lsm

 \begin{tabbing}
 \hspace{5mm} \= \hspace{5mm} \= \hspace{5mm} \= \hspace{60mm} \= \kill

 \textit{Factor the near-field.}\\
 \textbf{loop} over levels $l$ in the binary tree $\mathcal{T}_{\rm near}$\\
\> \textbf{loop} over boxes $\tau$ on level $l$\\
\> \> \textbf{if} $\tau$ is a leaf box,\\
\> \> \> $I_{\tau}$ is the indices of points in $\tau$.\\
\> \> \> Compute the low rank factorization of $\mtx{A}_{op}(I_{\tau},:)$ via Algorithm \ref{alg:geo}.\\
\> \> \> The interpolation matrix $\mtx{P}_\tau$ and index vector $J_\tau$ are returned.\\
\> \> \> $\mtx{U}_\tau = \mtx{P}_\tau$. \\
\> \> \textbf{else} \\
\> \> \> Let $\sigma_1$ and $\sigma_2$ denote the children of $\tau$.\\
\> \> \> $I_{\tau} = [J_{\sigma_1},J_{\sigma_2}]$. \textit{($J_{\sigma_j}$ denotes the index vector from factoring $\sigma_j$)}\\
\> \> \> Compute the low rank factorization of $\mtx{A}_op(I_{\tau},:)$ via Algorithm \ref{alg:geo}.\\
\> \> \> The interpolation matrix $\mtx{P}_\tau$ and index vector $J_\tau$ are returned.\\
%
\> \> \> $\mtx{U}_{\tau} = \mtwo{\mtx{U}_{\sigma_1}}{\mtx{0}}{\mtx{0}}{\mtx{U}_{\sigma_2}}\mtx{P}_\tau$.\\
\> \> \textbf{end if}\\
\> \textbf{end loop}\\
 \textbf{end loop}\\ 
$\mtx{L}_{\rm near} = \mtx{U}_1$, $\mtx{R}_{\rm near} = \mtx{A}_{op}(J_{\rm near},:)$, and $J_{\rm near} = J_1$\\

 \\
 
\textit{Extract factors for far-field from HBS representation of $\mtx{A}_{oo}$.}\\
\textbf{loop} over boxes $\tau$ from smallest to largest \\
\> \textbf{if} $\tau\subset \Gamma_{\rm far}$, \\
\>\> \textbf{if} $\tau$ is a leaf box \\
\> \> \> Set $\mtx{U}_\tau = \mtx{P}_\tau$ \textit{(the interpolation matrix 
      from the HBS representation of $\mtx{A}_{oo}$)}.\\
     \> \> \textbf{else} \\
\> \> \> Let $\sigma_1$ and $\sigma_2$ denote the children of $\tau$.\\
\> \> \> $\mtx{U}_{\tau} = \mtwo{\mtx{U}_{\sigma_1}}{\mtx{0}}{\mtx{0}}{\mtx{U}_{\sigma_2}}\mtx{P}_\tau$.\\
\> \> \textbf{end if}\\
\> \textbf{end if}\\
 \textbf{end loop}\\
 Let $T=\{\tau_1,\ldots,\tau_m\}$ denote the set containing the largest boxes touched in 
 the previous loop.\\
 Let $\mtx{L}_{\rm far}$ denote the block diagonal matrix with the $\mtx{U}_\tau$ matrices.\\
 \\

    Use algorithm \ref{alg:recompress} to recompress the far-field factors. Then
    $\mtx{L}_{\rm far} = \mtx{L}_{\rm far}\mtx{L}$ and $\mtx{R}_{\rm far} = \mtx{A}_{op}(J,:)$.

\\ 

   \textit{Do one more compression to remove extra degrees of freedom}\\
   $\mtx{L} = \mtwo{\mtx{L}_{\rm far}}{0}{0}{\mtx{L}_{\rm near}}$, 
   $\mtx{R} = \vtwo{\mtx{R}_{\rm far}}{\mtx{R}_{\rm near}}$, and 
   $J_{\rm tot} = [J_{\rm far}, J_{\rm near}]$.\\ 
   $[\mtx{P},J] = id(R);$\\

   \\
   
   $\mtx{L}_{op} = \mtx{L}\mtx{P}$, \qquad 
   $\mtx{R}_{op} = \mtx{A}_{op}(J,:)$.
   
 \end{tabbing}
\end{minipage}}
\end{center}
\end{figure}


In the case where $N_c<<N_k$, the total compression cost of $\mtx{A}_{op}$ 
and $\mtx{A}_{kc}$ can be further reduced by combining the two 
and only performing one transversal of the binary tree on $\Gamma_o$.

\subsection{Computational cost for the precomputation}
\label{sec:precomp_calc}

Recall that the proposed direct solver is comprised of two steps: precomputation and solve. 
The precomputation step is more expensive than the solve step but it needs to only be computed
once.  As stated in the beginning of this section, the precomputation consists of computing 
the low rank factorization of the update matrix $\mtx{Q}$, $\mtx{A}^{-1}\mtx{L}$ and 
inverting $\left(\mtx{I}+\mtx{R}\mtx{A}^{-1}\mtx{L}\right)$.  This section details the 
computational cost of the precomputation step. 

The computational cost of constructing the low rank factorization of 
$\mtx{Q}$ is $O((N_o+N_p)k_Q)$, where $k_Q = \max \{k_{\rm HBS},k_p\}$, $k_{\rm HBS}$ 
denotes the rank of the HBS factors and $k_p$ denotes the rank of the 
$\mtx{A}_{op}$ factorization.  The value of $k_p$ depends on the relationship
of $\Gamma_p$ and $\Gamma_o$.  

Let the low rank factors $\mtx{L}$ and $\mtx{R}^T$ of $\mtx{Q}$ have size
$N \times k$ where $N = N_o+N_p$ and $k = k_{kc}+N_c+k_{pk}+k_{op}$.  Constructing $\mtx{A}^{-1}\mtx{L}$ can 
be done in a block fashion by 

\begin{align*}
 \mtx{A}^{-1}\mtx{L} &= \mtwo{\mtx{A}_{oo}^{-1}}{0}{0}{\mtx{A}_{pp}^{-1}} \begin{bmatrix}\begin{pmatrix} -\mathbf{L}_{kc} &\mathbf{0}\\ 0& -\mathbf{B}_{cc} \end{pmatrix} & \mathbf{L}_{op}\\
\mathbf{L}_{pk} & \mathbf{0} \\
\end{bmatrix} \\
&=\begin{bmatrix}
-\mbf{A}_{oo}^{-1}\begin{pmatrix} \mathbf{L}_{kc} &\mathbf{0}\\ 0& \mathbf{B}_{cc} \end{pmatrix} & \mbf{A}_{oo}^{-1} \mathbf{L}_{op}\\
\mbf{A}_{pp}^{-1}\mbf{L}_{pk} & \mbf{0}
\end{bmatrix}.\end{align*}

The computational cost of computing the upper left and right blocks is $O\left(N_o(k_{kc}+N_c)\right)$ 
and $O(N_ok_{op})$, respectively. If $N_p$ is small, it is efficient to compute 
$\mtx{A}_{pp}^{-1}$ via dense linear algebra for $O(N_p^3)$ cost.  For $N_p$ large,
an approximate inverse of $\mtx{A}_{pp}$ can be constructed via fast direct solver for 
$O(N_p)$ computational cost.  The computational cost of evaluating $\mtx{A}_{pp}^{-1}\mtx{L}_{pk}$
is $O(N_p^2k_{pk})$ via dense linear algebra and $O(N_pk_{pk})$ via fast linear algebra.
Thus for $k_{op}$ small and $N_c$ constant, the cost of constructing $\mtx{A}^{-1}\mtx{L}$ is 
linear with respect to $N_o$.  

For problems of interest, the matrix $\left(\mtx{I}+\mtx{R}\mtx{A}^{-1}\mtx{L}\right)$ 
is small enough to be inverted rapidly via dense linear algebra for $O(k^3)$ computational 
cost. 

%

\subsection{Computational cost of the solve step}
\label{sec:solve_calc}

The solve step consists of applying the approximate inverse of $\mtx{A}$ to the 
the vector $\vct{f}_{\rm ext}$ in equation (\ref{eq:extended_system}),
matrix vector multiplications and vector addition. 
As discussed in section \ref{sec:precomp_calc},
$\mtx{A}^{-1}$ is applied by blocks for a cost $O(N_o+ N_p^2)$.  
Applying $\mtx{R}$, $\left(\mtx{I}+\mtx{R}\mtx{A}^{-1}\mtx{L}\right)^{-1}$,
and $\mtx{A}^{-1}\mtx{L}$ have computational cost $O((N_o+N_p)k)$, $O(k^2)$,
and $O((N_o+N_p)k)$. Thus the total computation cost of the 
solve step is $O((N_o+N_p)k+k^2+ N_o+N_p^2)$.  Thus the solve step
is linear with respect to $N_o$.

\section{Numerical experiments}
\label{sec:numerics}

\begin{figure}[H]
\begin{center}
\begin{picture}(100,100)(100,50)
\put(20,20){\includegraphics[scale=0.35]{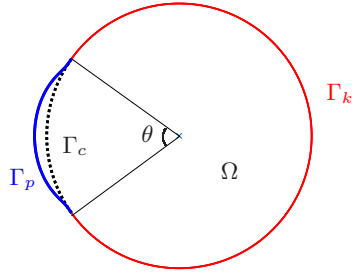}}
\put(140,80){\footnotesize${\Omega}$}
\put(180,110){\footnotesize\color{red}${\Gamma}_k$}
\put(80,90){\footnotesize${\Gamma}_c$}
\put(60,78){\footnotesize\color{blue}${\Gamma}_p$}
\put(110,94){\footnotesize$\theta$}
\qbezier(119,93)(117,96)(119,99)
\end{picture}
\end{center}
\caption{The circle with a bump of central angle $\theta$ geometry.}
\label{fig:test3_geom}
\end{figure}

This section illustrates the performance of the new fast direct solver for 
three types of locally-perturbed geometries.  The geometries under consideration
are \\
\noindent
\textbf{Square with a nose:} The original boundary $\Gamma_o$ is a square
with corners rounded via the method in \cite{doi:10.1137/15M1028248}.  The local
perturbation $\Gamma_p$ is a rectangle with rounded corners that is attached to $\Gamma_o$
again using the method in \cite{doi:10.1137/15M1028248} to smooth. 
Figure \ref{fig:test1_geom_smooth} illustrates the geometry.  The length of 
the local perturbation rectangle is fixed but the height $d$ can vary depending 
on experiment.  Composite Gaussian quadrature is used to discretize the integral equation.\\
\noindent 
\textbf{Circle with a bump:}  The original geometry is a 
circle. The perturbed geometry replaces an arc of the circle by a smooth bump with central angle $\theta$.
Figure \ref{fig:test2_geom} illustrates the geometry.  This geometry is considered in \cite{2016_update}. 
Trapezoidal rule is used to discretize the integral equation.\\
\noindent
\textbf{Star with refined panels:} The original geometry is a star geometry discretized with composite
Gaussian quadrature (see Figure \ref{fig:test3_geom}(a)).  To create the perturbed boundary, three Gaussian 
panels illustrated in Figure \ref{fig:test3_geom}(b) are replaced with more panels. Figure \ref{fig:test3_geom}(b) 
illustrates a local perturbation $\Gamma_p$ consisting of six panels. \\

To ensure the solution technique preserves accuracy, we test it on 
problems with a known exact solution given by 
$u_{exact}(\vec{x})=\sum_{j=1}^{10}q_j G(\vec{x},\vec{s}_j)$, 
where $\{\vec{s}_j\}_{j=1}^{10}$ are point charges placed 
outside of $\Omega$ and $\{q_j\}_{j=1}^{10}$ are the charge values.  
We define the the relative error to be
\begin{equation}
\label{equ:num_res_rel_err}
E = \frac{\|\mbf{u}_{exact}-\mbf{u}_{new}\|_2}{\|\mbf{u}_{exact}\|_2},
\end{equation}
where the vector $\vct{u}_{new}$ and $\vct{u}_{exact}$ contain the approximate and exact solution, respectively, at ten points 
$\{\vec{t}\}_{j=1}^{10}\in\Omega$.  For all problems under consideration, the geometries 
are fully resolved.  Thus, with the tolerance of the 
compression schemes set to $\epsilon  = 10^{-10}$, the 
relative error $E$ is approximately $10^{-9}$ for all 
choices of $N_o$ and $N_p$.  

For each geometry, we report the following: \\
\noindent
$N_o$ the number of discretization points on the original geometry.\\
\noindent
 $N_p$ the number of discretization points on the added geometry.\\
\noindent
 $T_{new,\,p}$ the time in seconds for the precomputation 
 step for the new solver. \\
\noindent
 $T_{hbs,\,p}$ the time in seconds to construct a new HBS solver.\\
\noindent
 $r_p:= \frac{T_{new,\,p}}{T_{hbs,\,p}}$.\\
\noindent
  $T_{new,\,s}$ the time in seconds to apply the new 
  solver.  \\
\noindent
  $T_{hbs,\,s}$ the time in seconds to apply the HBS approximation of 
  the inverse.\\
\noindent
 $r_s:=\frac{T_{new,\,s}}{T_{hbs,\,s}}$. \\
 
The ratios $r_p$ and $r_s$ illustrate the performance of the new solver relative 
to building a new HBS solver from scratch for the boundary value problem
on the perturbed geometry.

\subsection{Square with a nose}
\label{sec:nose}

\begin{figure}[H]
\begin{center}
\begin{picture}(150,100)(100,50)
\put(20,20){\includegraphics[scale=0.35]{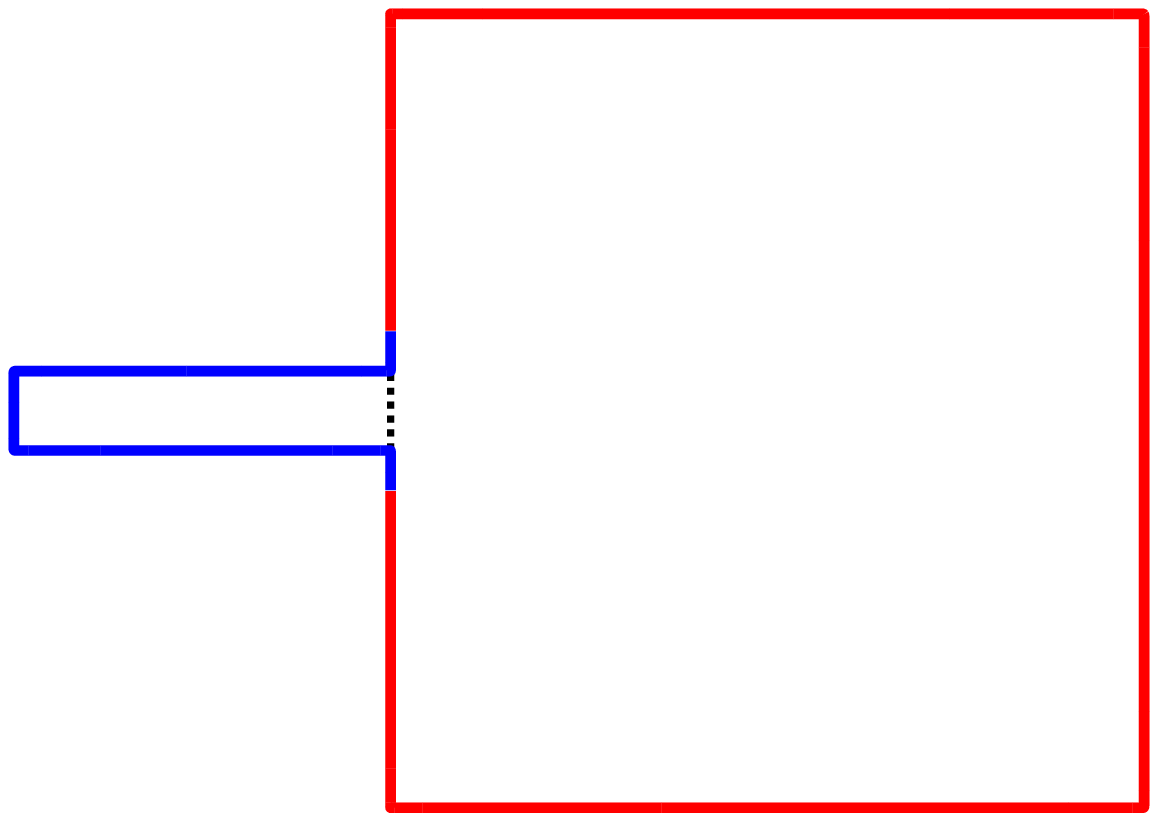}}
\put(219,70){\includegraphics[scale=0.11]{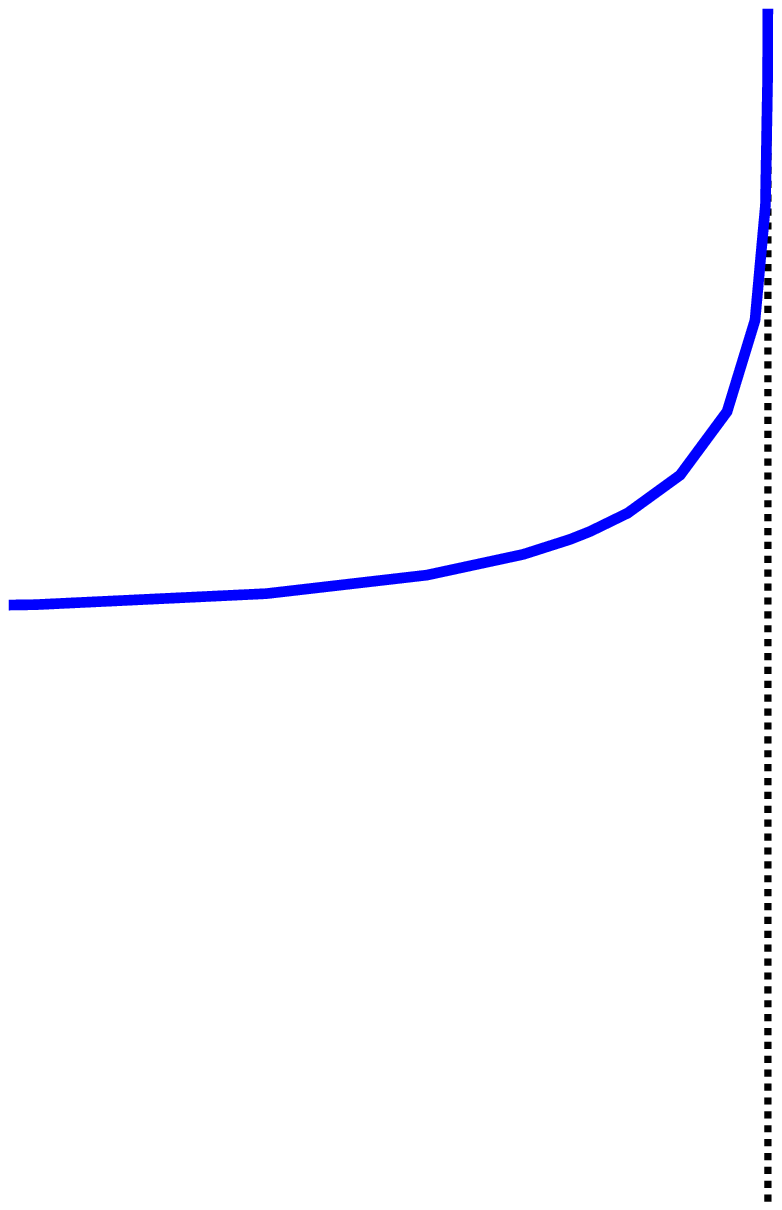}}
\put(140,72){\footnotesize${\Omega}$}
\put(116,122){\footnotesize\color{red}${\Gamma}_k$}
\put(107,90){\footnotesize${\Gamma}_c$}
\put(240,80){\tiny${\Gamma}_c$}
\put(235,99){\tiny\color{blue}${\Gamma}_p$}
\put(68,78){\footnotesize\color{blue}${\Gamma}_p$}
\put(54,94){\footnotesize$d$$\{$}
\put(102,100){\circle{6}}
\qbezier(102,100)(155,120)(220,100)
\put(220,100){\vector(1,-1){2}}
\put(245,94){\circle{62}}
\end{picture}
\end{center}
\caption{
The square geometry with a smoothly
attached nose of height $d$.}
\label{fig:test1_geom_smooth}
\end{figure}

This section reports the performance of the new solver on the square with nose geometry illustrated in Figure \ref{fig:test1_geom_smooth}.
Two choices of ``nose'' height $d$ are considered: thinning and fixed.  For the 
thinning nose geometry, the height $d$ is decreased as $N_o$ increases allowing $N_c$ to remain 
constant.  For the fixed nose geometry, $N_c$ is thus increasing at the same rate as $N_o$. 
For both geometries, $N_p$ varies between 700 and 900 due to the corner-smoothing procedure.

\begin{figure}[H]
\begin{center}
\begin{picture}(300,190)(70,01)
\put(-20,06){\includegraphics[scale=0.4]{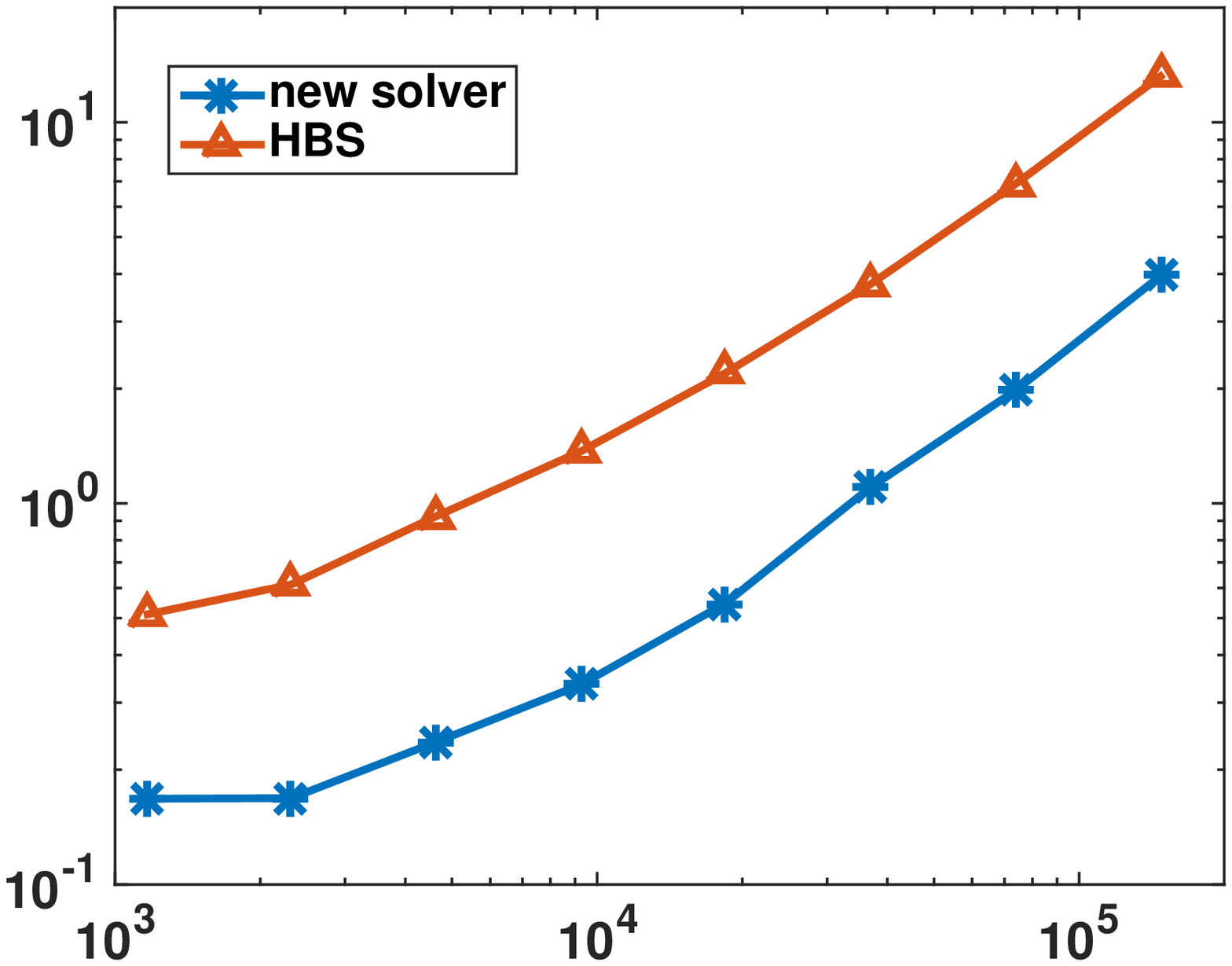}}
\put(200,06){\includegraphics[scale=0.4]{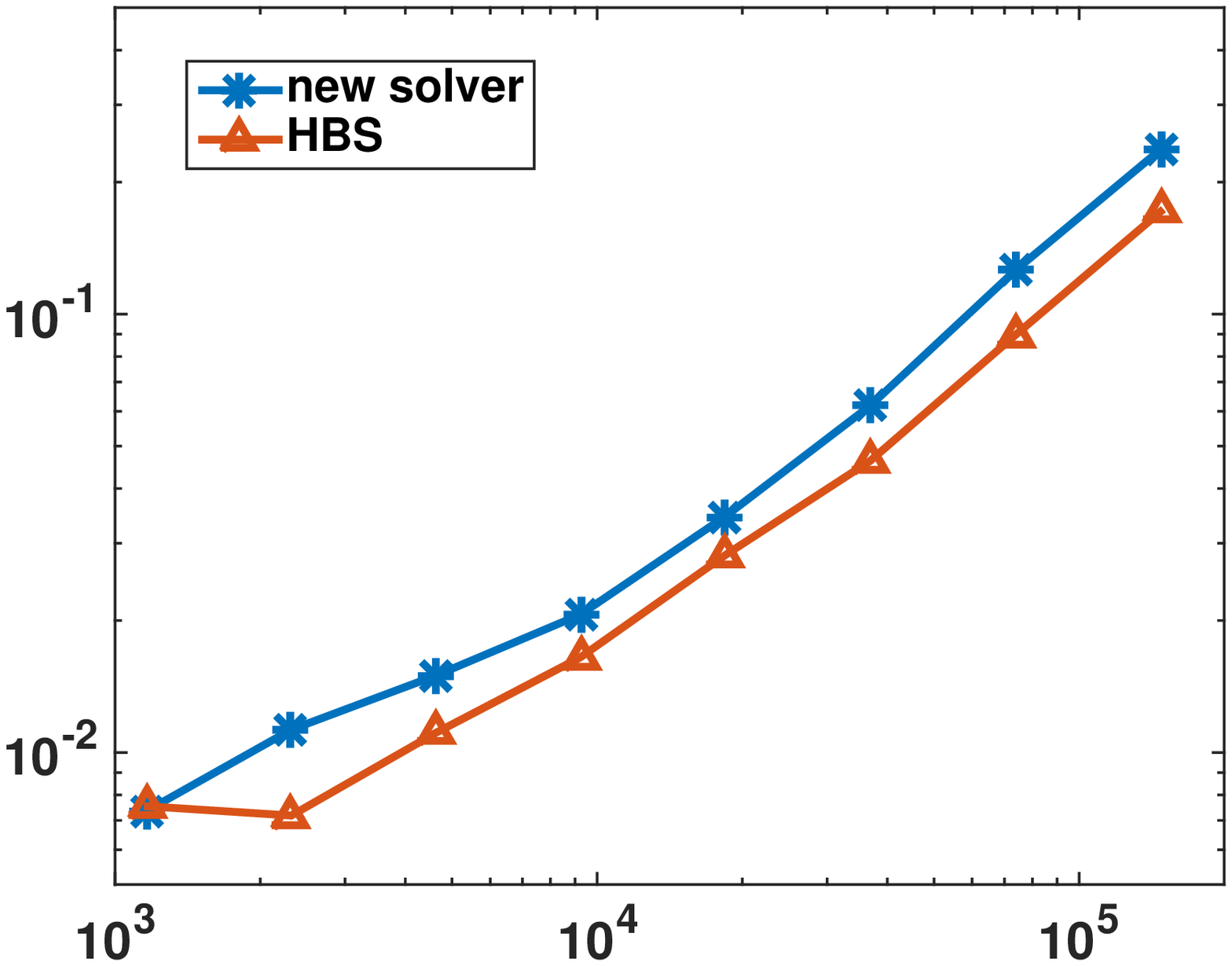}}
\put(50,175){ Precomputation}
\put(290,175){ Solve}
\put(-30,65){\rotatebox{90}{Time (sec)}}
\put(90,-10){$(a)$}
\put(300,-10){$(b)$}
\put(90,-0){$N_o$}
\end{picture}
\end{center}
\caption{A log-log plot of the time in seconds versus $N_o$ for the 
(a) precomputation and (b) solve steps of both the new and HBS solver
for the square with thinning nose geometry in Section \ref{sec:nose}.}
\label{fig:test1_shrink}
\end{figure}

 Figure \ref{fig:test1_shrink} presents log-log plots of the time in
 seconds versus $N_o$ for the (a) precomputation and (b) solve steps for
 the new solver and HBS solver for the thinning nose geometry. This figure 
 and the timings reported in Table \ref{tab:test1_shrink} illustrate that 
 the new solver does scale linearly with respect to $N_o$ for this problem.
 The entry $r_p$ in Table \ref{tab:test1_shrink} reports that the precomputation step of the new solver
 is approximately three times faster than the precomputation of the HBS solver.
 The solve step of the new solver is slower than the solve step
 of the HBS solver.  However given the much larger constant associated with
 the precomputation step, it would take 100 to 260 solves
 to make the new solver slower than building a new HBS solver from 
 scratch.

\begin{table}[H]
\begin{tabular}{l|c|c|c|c|c|c} 
\hline
$N_o$  & $\,T_{new,\,p}\,$& $\,T_{hbs,\,p}\,$ & $\,r_p\,$ & $\,T_{new,\,s}$\,& \,$T_{hbs,\,s}$ \,& $\,r_s$ \\
  \hline
1168      &  1.68e-01 &  5.12e-01 & 3.28e-01 &  7.38e-03 &  7.54e-03 & 9.79e-01 \\ 
  \hline
2320      &  1.69e-01 &  6.13e-01 & 2.75e-01 &  1.13e-02 &  7.18e-03 & 1.57e+00 \\  
  \hline
4624      &  2.36e-01 &  9.24e-01 & 2.56e-01 &  1.50e-02 &  1.11e-02 & 1.35e+00 \\ 
  \hline
9232      &  3.34e-01 &  1.37e+00 & 2.44e-01 &  2.06e-02 &  1.65e-02 & 1.25e+00 \\
  \hline
18448      &  5.47e-01 &  2.20e+00 & 2.49e-01 &  3.46e-02 &  2.82e-02 & 1.23e+00 \\ 
  \hline
36880      &  1.10e+00 &  3.76e+00 & 2.93e-01 &  6.18e-02 &  4.63e-02 & 1.34e+00 \\  
  \hline
73744      &  1.98e+00 &  6.88e+00 & 2.87e-01 &  1.26e-01 &  8.96e-02 & 1.41e+00 \\
  \hline
147472      &  3.95e+00 &  1.32e+01 & 2.99e-01 &  2.37e-01 &  1.71e-01 & 1.39e+00 \\
  \hline

   \end{tabular} 
   
   \caption{Timing results for the square with thinning nose geometry in Section \ref{sec:nose}.}
   \label{tab:test1_shrink}
  
\end{table}

\begin{figure}[H]
\begin{center}
\begin{picture}(300,190)(70,01)
\put(-20,06){\includegraphics[scale=0.4]{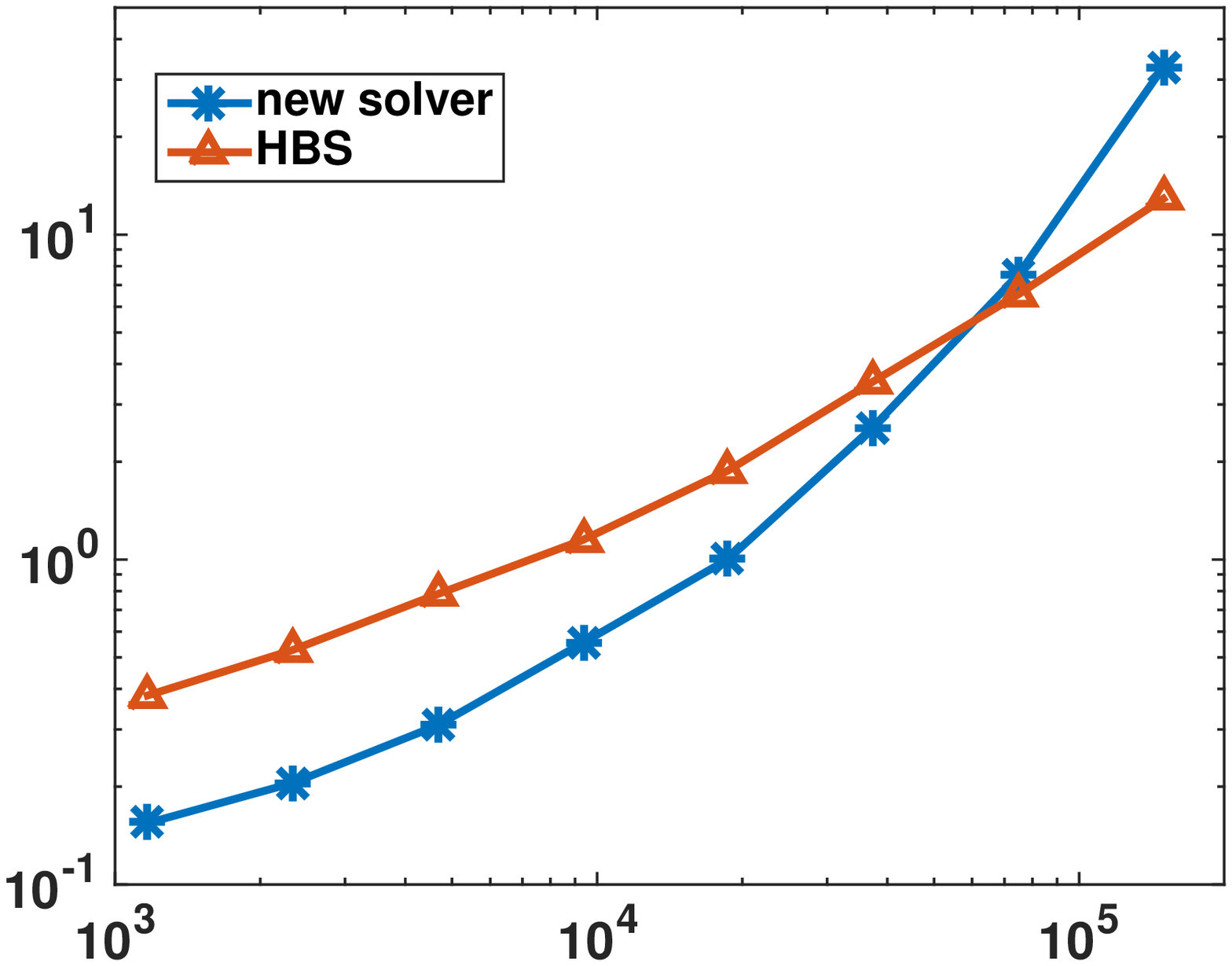}}
\put(200,06){\includegraphics[scale=0.4]{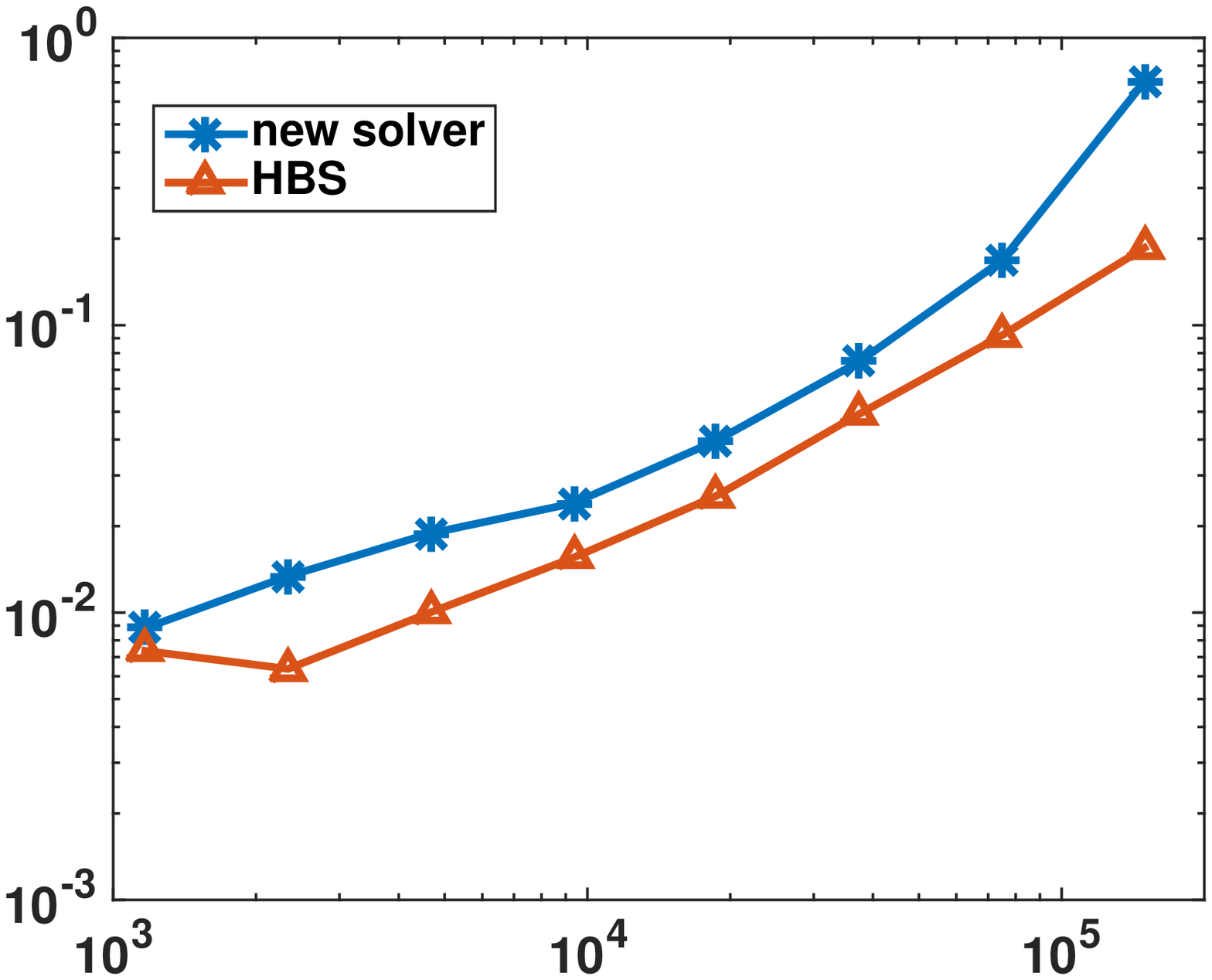}}
\put(50,175){ Precomputation}
\put(290,175){ Solve}
\put(-30,65){\rotatebox{90}{Time (sec)}}
\put(90,-10){$(a)$}
\put(300,-10){$(b)$}
\put(90,1){$N_o$}
\end{picture}
\end{center}
\caption{A log-log plot of the time in seconds versus $N_o$ for the 
(a) precomputation and (b) solve steps of both the new and HBS solver
for the square with fixed nose geometry in Section 
  \ref{sec:nose}.}
\label{fig:test1_noshrink}

\end{figure}

 Figure \ref{fig:test1_noshrink} presents log-log plots of the time in
 seconds versus $N_o$ for the (a) precomputation and (b) solve steps for
 the new solver and HBS solver for the fixed nose geometry. As predicted in 
 section \ref{sec:precomp_calc}, the computational cost for both steps 
 does not scale linearly with respect to $N_o$.  
Table \ref{tab:test1_no_shrink} reinforces this statement.
Specifically, notice that the time to build the HBS solver ($T_{hbs,p}$)
scales linearly with $N_o$ while the time for the new solver
does not, and thus $r_p$ does not approach a constant.  The 
lack of linear scaling is a direct result of the fact that 
the rank of the update matrix $\mbf{Q}$ is dependent on $N_o$.

\begin{table}[H]
\begin{tabular}{l|c|c|c|c|c|c} 
\hline
$N_o$  & $\,T_{new,\,p}\,$& $\,T_{hbs,\,p}\,$ & $\,r_p\,$ & $\,T_{new,\,s}$\,& \,$T_{hbs,\,s}$ \,& $\,r_s$ \\
 \hline

1168      &  1.55e-01 &  3.82e-01 & 4.05e-01 &  8.82e-03 &  7.36e-03 & 1.20e+00 \\ 
  \hline
2336      &  2.05e-01 &  5.27e-01 & 3.89e-01 &  1.34e-02 &  6.37e-03 & 2.10e+00 \\ 
  \hline
4672      &  3.09e-01 &  7.83e-01 & 3.95e-01 &  1.89e-02 &  1.00e-02 & 1.89e+00 \\ 
  \hline
9344      &  5.54e-01 &  1.15e+00 & 4.82e-01 &  2.40e-02 &  1.55e-02 & 1.55e+00 \\ 
  \hline
18688      &  1.00e+00 &  1.88e+00 & 5.34e-01 &  3.96e-02 &  2.54e-02 & 1.56e+00 \\  
  \hline
37376      &  2.56e+00 &  3.53e+00 & 7.24e-01 &  7.46e-02 &  4.90e-02 & 1.52e+00 \\
  \hline
74752      &  7.51e+00 &  6.55e+00 & 1.15e+00 &  1.68e-01 &  9.20e-02 & 1.82e+00 \\ 
  \hline
149504      &  3.277e+01 &  1.290e+01 & 2.541e+00 &  7.038e-01 &  1.858e-01 & 3.787e+00 \\
  \hline
   \end{tabular} 
   \caption{Timing results for the square with fixed nose geometry in Section 
  \ref{sec:nose}.}
   \label{tab:test1_no_shrink}
\end{table}

\subsection{Circle with a bump}
\label{sec:circlebump}
\begin{figure}[H]
\begin{center}
\begin{picture}(300,100)(01,01)
\put(-100,-20){\includegraphics[scale=0.35]{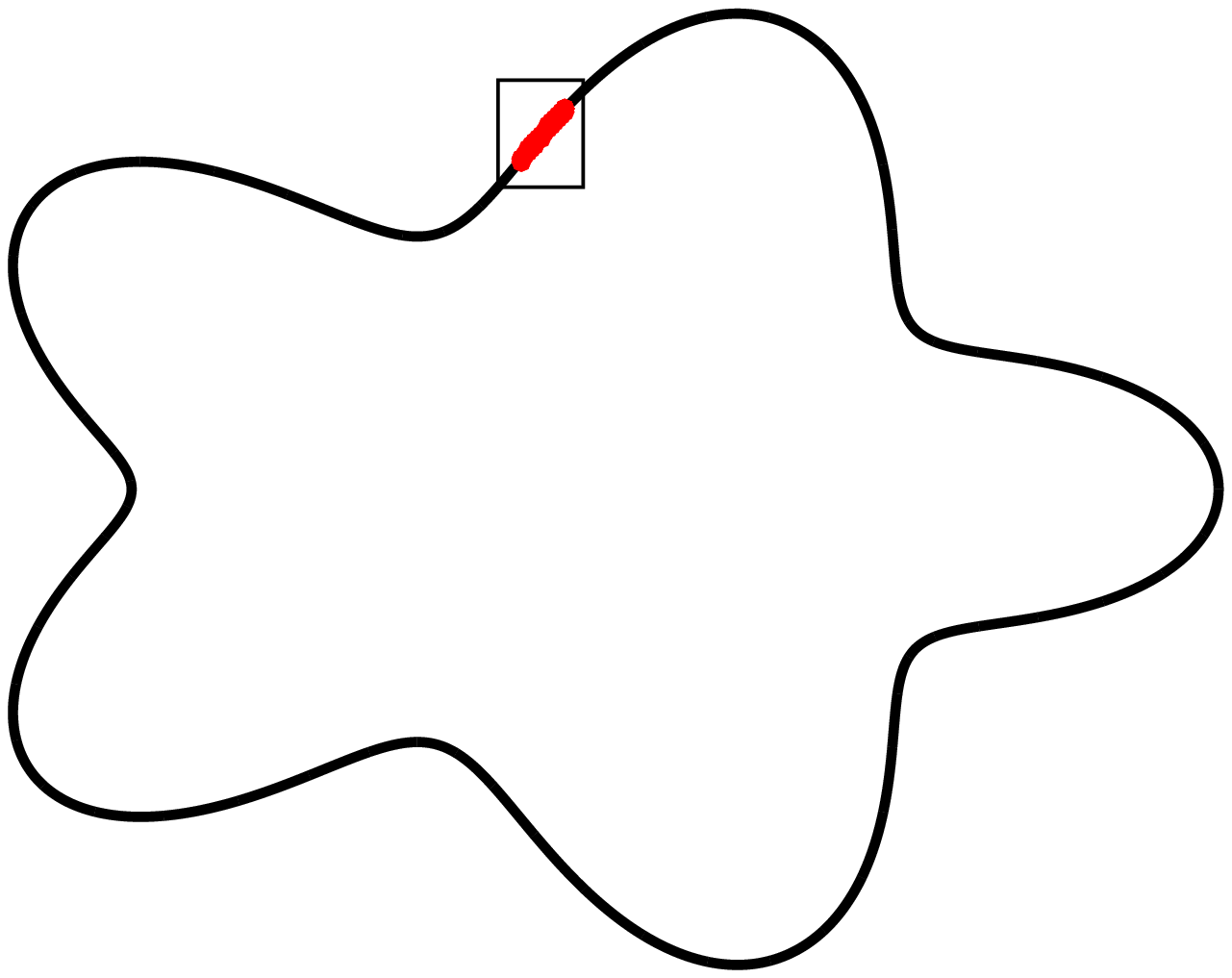}}
\put(40,20){$\Gamma_k$}
\put(-10,-10){(a)}
\put(80,01){\includegraphics[scale=0.25]{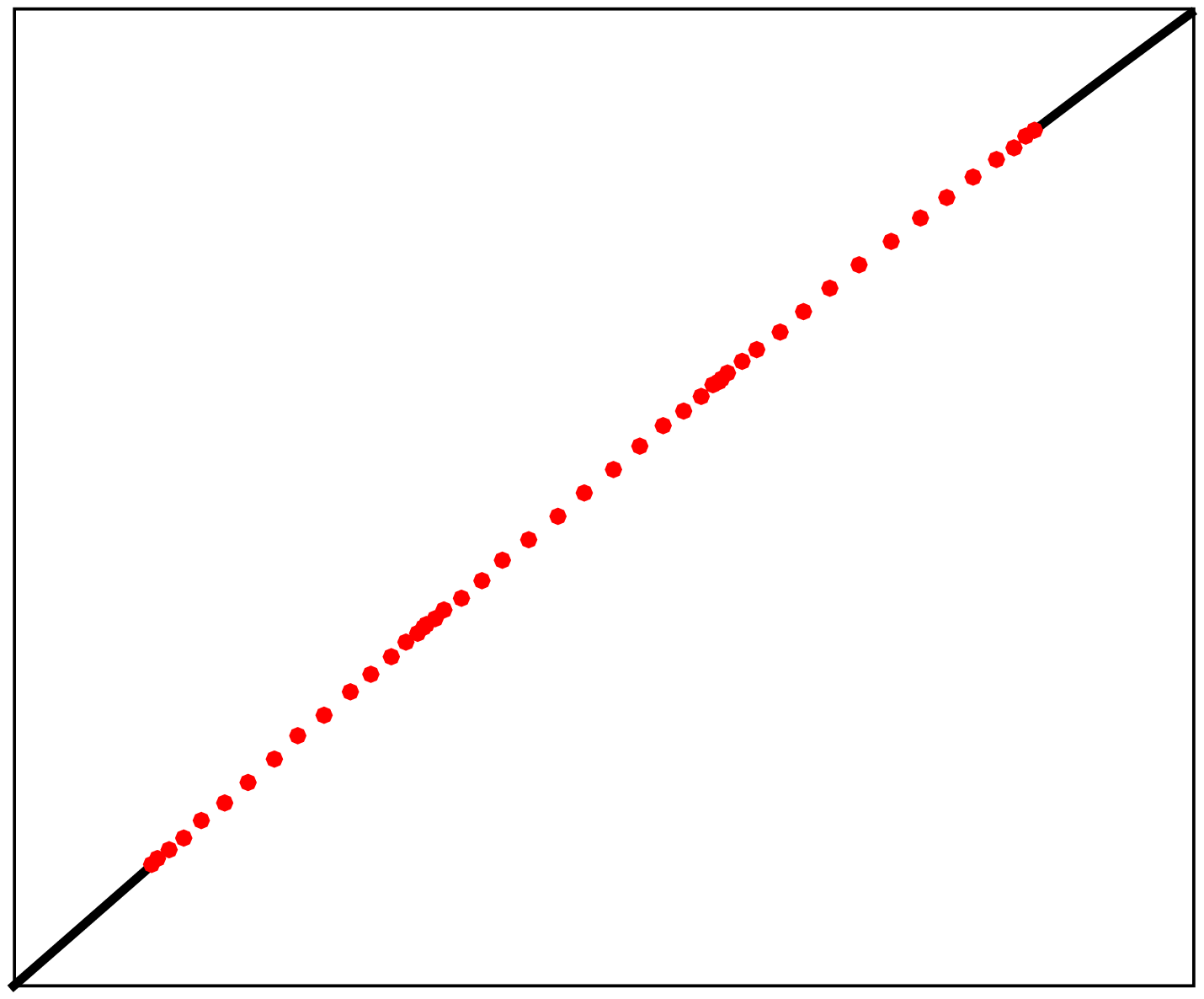}}
\put(140,-10){\footnotesize (b)}
\put(220,01){\includegraphics[scale=0.25]{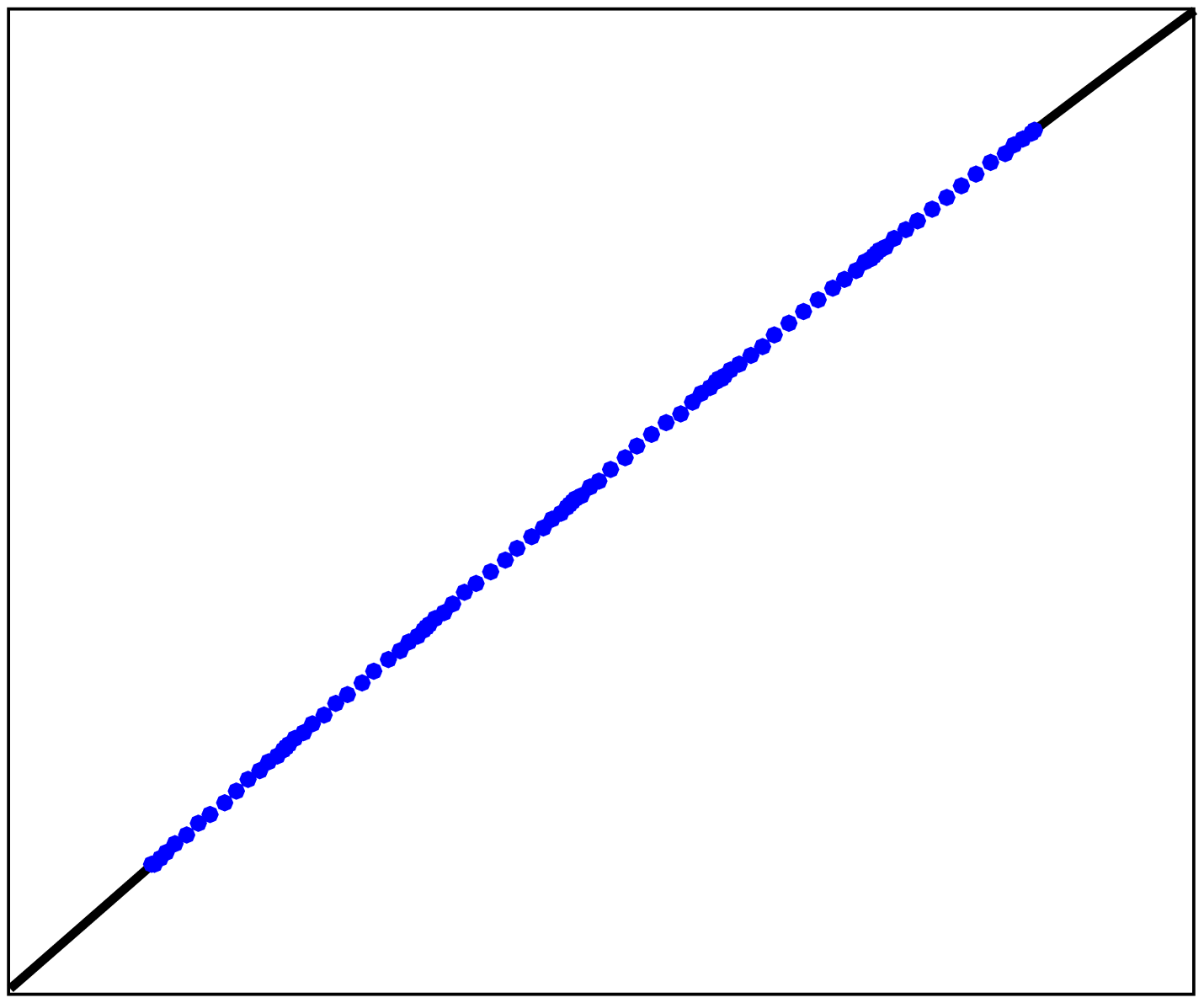}}
\put(280,-10){\footnotesize (c)}
\put(140,55){\color{red}$\Gamma_c$}
\put(280,55){\color{blue}$\Gamma_p$}

\end{picture}
\end{center}
\caption{(a) The star geometry with the portion of the boundary to 
be refined boxed.  (b) The three Gaussian panels in the boxed region
from the original discretization. (c) The six Gaussian panels that 
replaced the original three panels.
}
\label{fig:test2_geom}
\end{figure}

This section reports the performance of the new solver on the circle with bump geometry illustrated in Figure \ref{fig:test2_geom}.
Two choices of ``bump'' are considered: shrinking and fixed.  For the 
shrinking bump geometry, the angle $\theta$ decreases 
as $N_o$ increases so that $N_p=N_c=199$ independent of $N_o$.   For the fixed bump geometry, the 
angle $\theta$ remains fixed as $N_o$ increases.  Thus both $N_c$ and $N_p$ increase at the 
same rate as $N_o$.  These are the same geometries considered in \cite{2016_update}.

\begin{figure}[H]
\begin{center}
\begin{picture}(300,190)(70,01)
\put(-20,06){\includegraphics[scale=0.4]{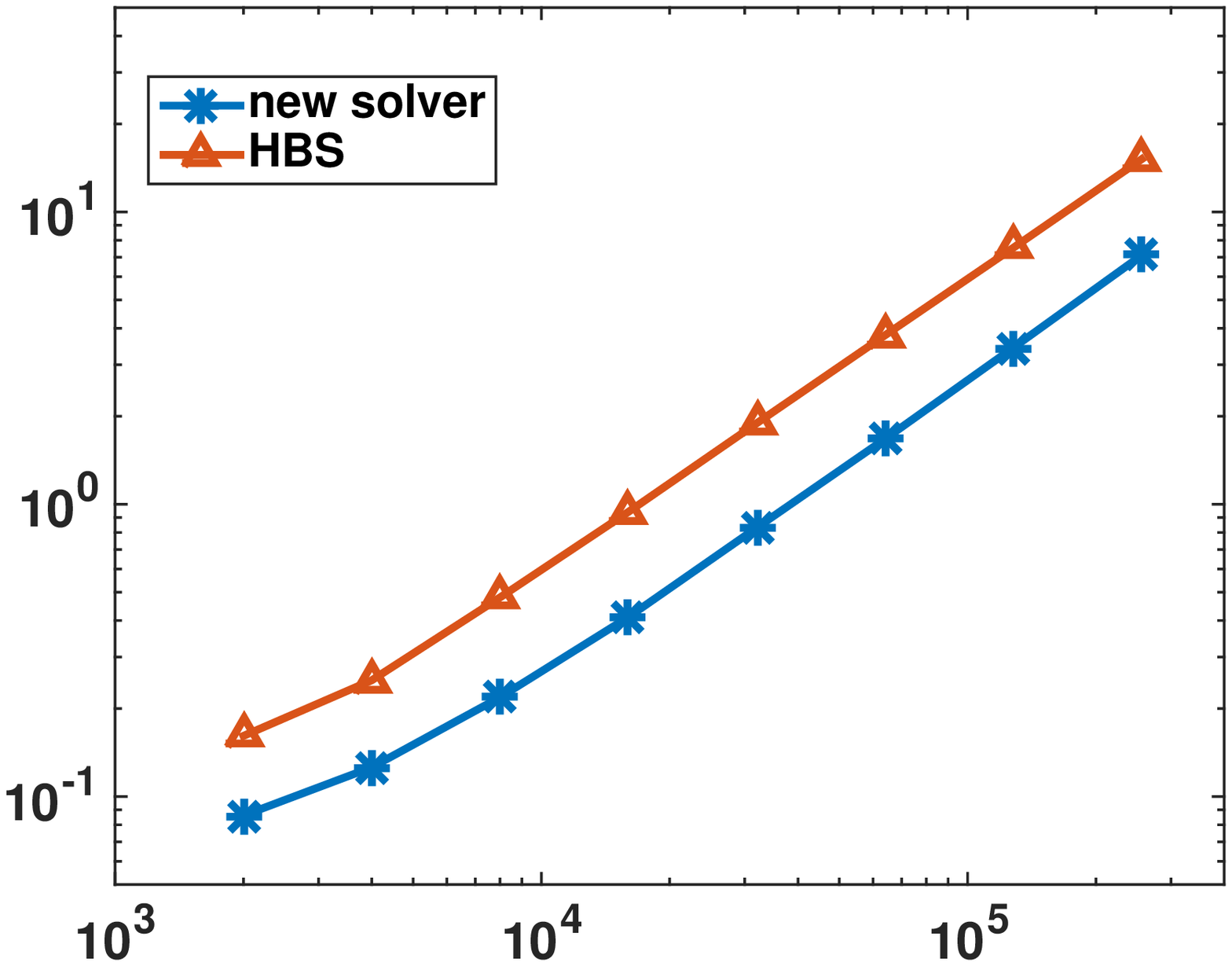}}
\put(200,06){\includegraphics[scale=0.4]{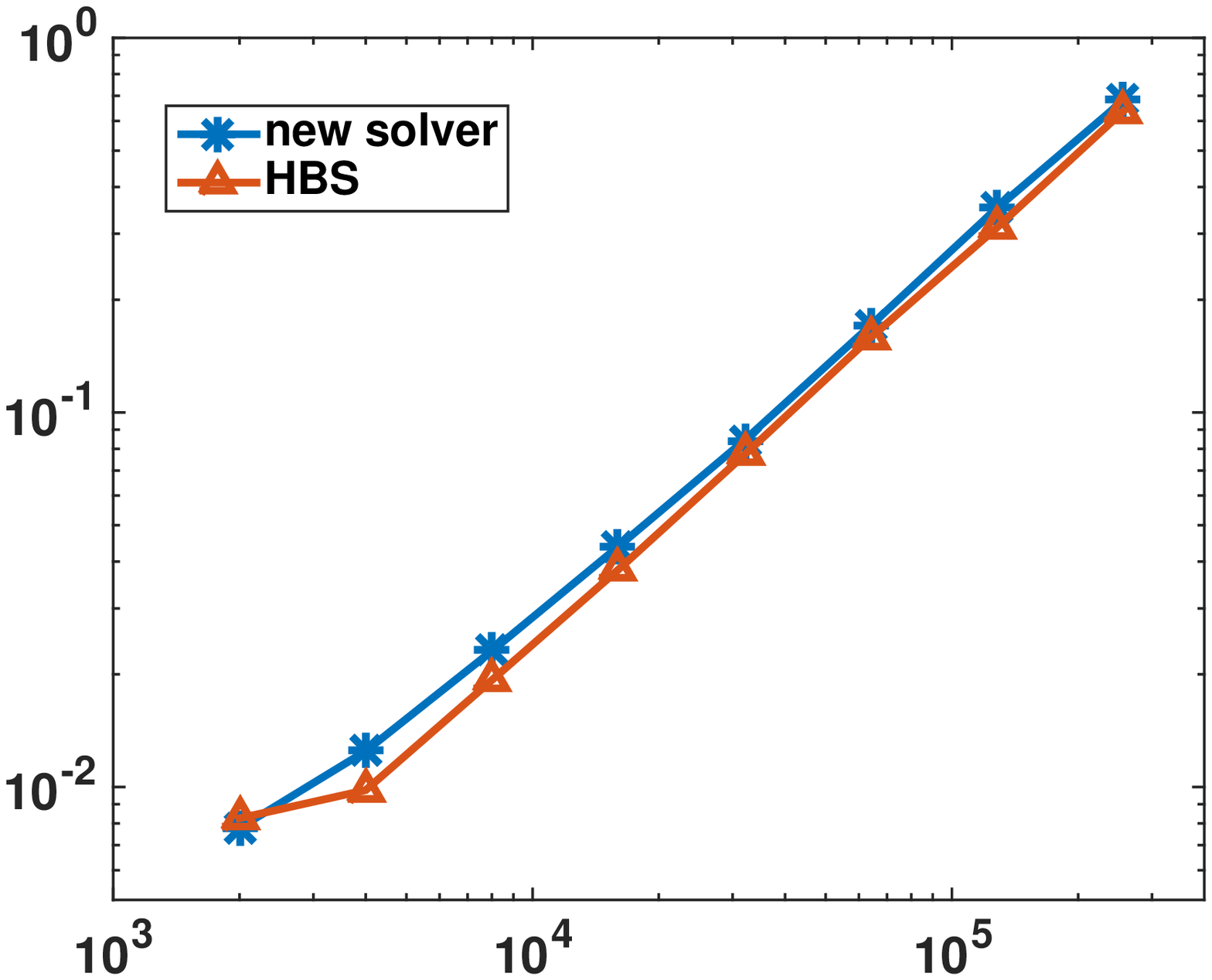}}
\put(50,175){ Precomputation}
\put(290,175){ Solve}
\put(-30,65){\rotatebox{90}{Time (sec)}}
\put(90,-10){$(a)$}
\put(300,-10){$(b)$}
\put(90,1){$N_o$}
\end{picture}
\end{center}
\caption{A loglog plot of the time in seconds versus $N_o$ for the 
(a) precomputation and (b) solve steps of both the new and HBS solver
for the circle with shrinking bump geometry in Section  \ref{sec:circlebump}. }
\label{fig:test3_shrink}

\end{figure}

 Figure \ref{fig:test3_shrink} presents log-log plots of the time in
 seconds versus $N_o$ for the (a) precomputation and (b) solve steps for
 the new solver and HBS solver for the shrinking bump geometry. As reported in Table 
 \ref{tab:test3_shrink}, the precomputation of the 
 new solver is linear and two times faster than the precomputation of the HBS solver.  Since $r_s \sim 1$, the cost of the 
solve step is nearly the same.  
Thus, the new solver is the more efficient choice for this geometry.

\begin{figure}[H]
\begin{center}
\begin{picture}(300,190)(70,01)
\put(-20,06){\includegraphics[scale=0.4]{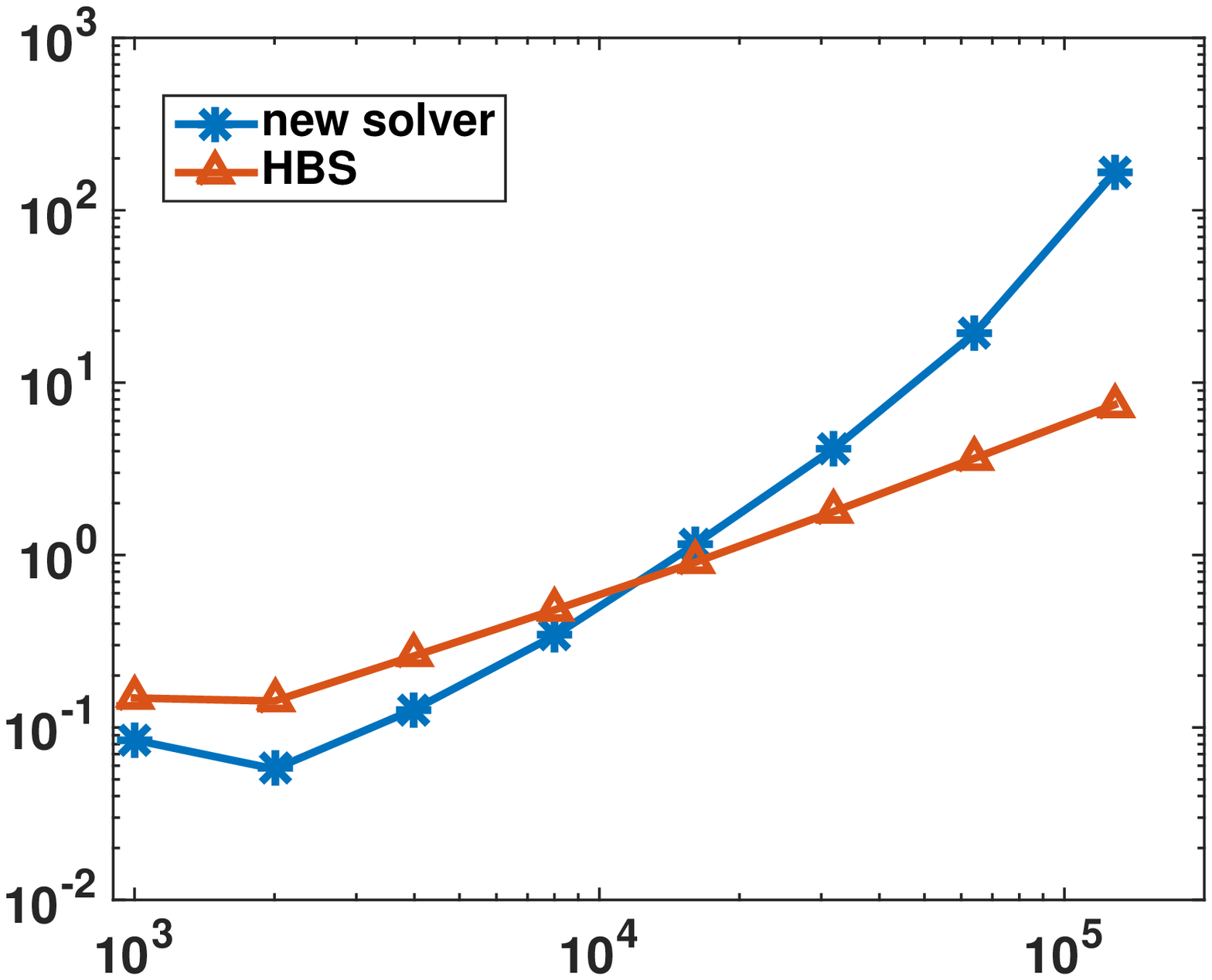}}
\put(200,06){\includegraphics[scale=0.4]{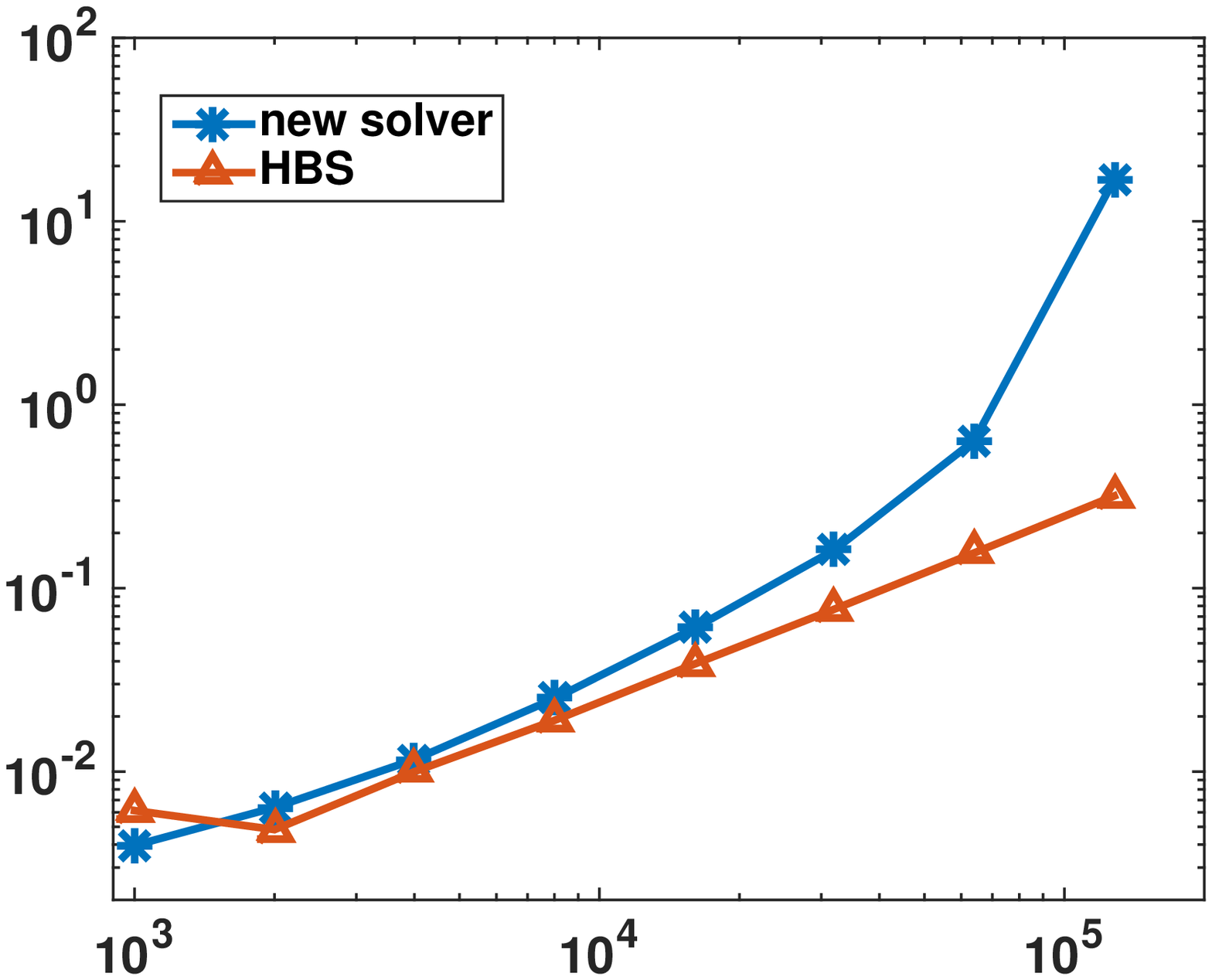}}
\put(50,175){ Precomputation}
\put(290,175){ Solve}
\put(-30,65){\rotatebox{90}{Time (sec)}}
\put(90,-10){$(a)$}
\put(300,-10){$(b)$}
\put(90,1){$N_o$}
\end{picture}
\end{center}
\caption{A log-log plot of the time in seconds versus $N_o$ for the 
(a) precomputation and (b) solve steps of both the new and HBS solver
for the circle with fixed bump geometry in Section \ref{sec:circlebump}. }
\label{fig:test3_noshrink}
\end{figure}

 Figure \ref{fig:test3_noshrink} presents log-log plots of the time in
 seconds versus $N_o$ for the (a) precomputation and (b) solve steps for
 the new solver and HBS solver for the fixed bump geometry.  For this geometry, linear scaling of the new solver is 
not expected since the size of 
the update matrix $\mbf{Q}$ grows with $N_o$ and $N_p$.  The timings reported in 
Table \ref{tab:test3_noshrink} and Figure \ref{fig:test3_shrink} support 
this statement.

\begin{table}[h]
\begin{tabular}{l|c|c|c|c|c|c} 
\hline
$N_o$  & $\,T_{new,\,p}\,$& $\,T_{hbs,\,p}\,$ & $\,r_p\,$ & $\,T_{new,\,s}$\,& \,$T_{hbs,\,s}$ \,& $\,r_s$ \\
 \hline
2000      &  8.61e-02 &  1.61e-01 & 5.34e-01&  7.77e-03 &  8.26e-03 & 9.42e-01 \\ 
  \hline
4000      &  1.26e-01 &  2.50e-01 & 5.04e-01 &  1.25e-02 &  9.81e-03 & 1.27e+00 \\
  \hline
8000      &  2.19e-01 &  4.80e-01 & 4.55e-01 &  2.31e-02 &  1.93e-02 & 1.20e+00 \\
  \hline
16000      &  4.11e-01 &  9.41e-01 & 4.37e-01 &  4.38e-02 &  3.80e-02 & 1.15e+00 \\
  \hline
32000      &  8.31e-01 &  1.89e+00 & 4.39e-01&  8.43e-02 &  7.71e-02 & 1.09e+00 \\ 
  \hline
64000      &  1.67e+00 &  3.78e+00 & 4.42e-01&  1.71e-01 &  1.59e-01 & 1.08e+00 \\ 
  \hline
128000      &  3.43e+00 &  7.54e+00 & 4.55e-01  &  3.51e-01 &  3.11e-01 & 1.13e+00 \\ 
  \hline
256000      &  7.10e+00 &  1.51e+01 & 4.69e-01 &  6.80e-01 &  6.36e-01 & 1.07e+00 \\ 
  \hline

   \end{tabular} 
   
   \caption{Timing results for the circle with shrinking bump geometry in Section \ref{sec:circlebump}.}
   \label{tab:test3_shrink}
  
\end{table}

\begin{table}[H]
\begin{tabular}{l|c|c|c|c|c|c|c} 
\hline
$N_o$  & $N_p$& $\,T_{new,\,p}\,$& $\,T_{hbs,\,p}\,$ & $\,r_p\,$ & $\,T_{new,\,s}$\,& \,$T_{hbs,\,s}$ \,& $\,r_s$ \\
\hline
1000      &49&  8.46e-02 &  1.48e-01 & 5.71e-01  &  3.96e-03 &  6.14e-03 & 6.46e-01 \\ 
  \hline
2000      &99&  5.74e-02 &  1.42e-01 & 4.04e-01  &  6.37e-03 &  4.83e-03 & 1.32e+00 \\ 
  \hline
4000      &199&  1.27e-01 &  2.60e-01 & 4.90e-01 &  1.16e-02 &  1.00e-02 & 1.16e+00 \\ 
  \hline
8000      &399&  3.40e-01 &  4.81e-01 & 7.08e-01 &  2.50e-02 &  1.90e-02 & 1.32e+00 \\ 
  \hline
16000      &799&  1.15e+00 &  9.05e-01 & 1.27e+00 &  6.04e-02 &  3.86e-02 & 1.57e+00 \\  
  \hline
32000      &1599&  4.16e+00 &  1.79e+00 & 2.32e+00 &  1.61e-01 &  7.68e-02 & 2.10e+00 \\ 
  \hline
64000      &3199&  1.93e+01 &  3.64e+00 & 5.31e+00 &  6.38e-01 &  1.56e-01 & 4.10e+00 \\
  \hline
128000    &6399  &  1.63e+02 &  7.43e+00 & 2.20e+01 &  1.69e+01 &  3.18e-01 & 5.31e+01 \\  
  \hline

   \end{tabular} 
   
   \caption{Timing results for the circle with fixed bump geometry in Section 
   \ref{sec:circlebump}.  }
   \label{tab:test3_noshrink}

\end{table}

\subsection{Star with refined panels}
\label{sec:star_panel}

\begin{table}[H]
\begin{tabular}{l|c|c|c|c|c|c|c} 
\hline
$N_p$  & $\frac{N_p}{N_o}$ & $\,T_{new,\,p}\,$& $\,T_{hbs,\,p}\,$ & $\,r_p\,$ & $\,T_{new,\,s}$\,& \,$T_{hbs,\,s}$ \,& $\,r_s$ \\
 \hline
\hline
96      &0.03  &1.34e-01 &  7.44e-01 & 1.80e-01 &  2.33e-02 &  2.49e-02 & 9.35e-01 \\  
  \hline
192      & 0.06 &  1.35e-01 &  6.65e-01 & 2.02e-01 &  2.34e-02 &  2.45e-02 & 9.55e-01 \\  
  \hline
384      & 0.12 &  1.67e-01 &  6.73e-01 & 2.48e-01 &  2.42e-02 &  2.48e-02 & 9.77e-01 \\ 
  \hline
768      & 0.24 &  2.34e-01 &  7.18e-01 & 3.25e-01 &  2.76e-02 &  2.55e-02 & 1.08e+00 \\ 
  \hline
1536      & 0.48  &  4.80e-01 &  8.33e-01 & 5.76e-01 &  5.87e-02 &  3.29e-02 & 1.78e+00 \\ 
  \hline
3072      & 0.96  &  1.29e+00 &  1.06e+00 & 1.22e+00 &  2.44e-01 &  3.94e-02 & 6.19e+00 \\ 
  \hline
6144      & 1.92 &  4.26e+00 &  1.63e+00 & 2.62e+00 &  1.70e+00 &  6.25e-02 & 2.72e+01 \\ 
  \hline
12288      & 3.84 &  2.79e+01 &  2.38e+00 & 1.17e+01 &  1.41e+01 &  9.81e-02 & 1.44e+02 \\ 
  \hline
     \end{tabular} 
   \caption{Timing results for the star with refined panels geometry in Section 
  \ref{sec:star_panel}. }
   \label{tab:test2}
\end{table}

\begin{figure}[H]
\begin{center}
\begin{picture}(300,190)(70,01)
\put(-20,06){\includegraphics[scale=0.4]{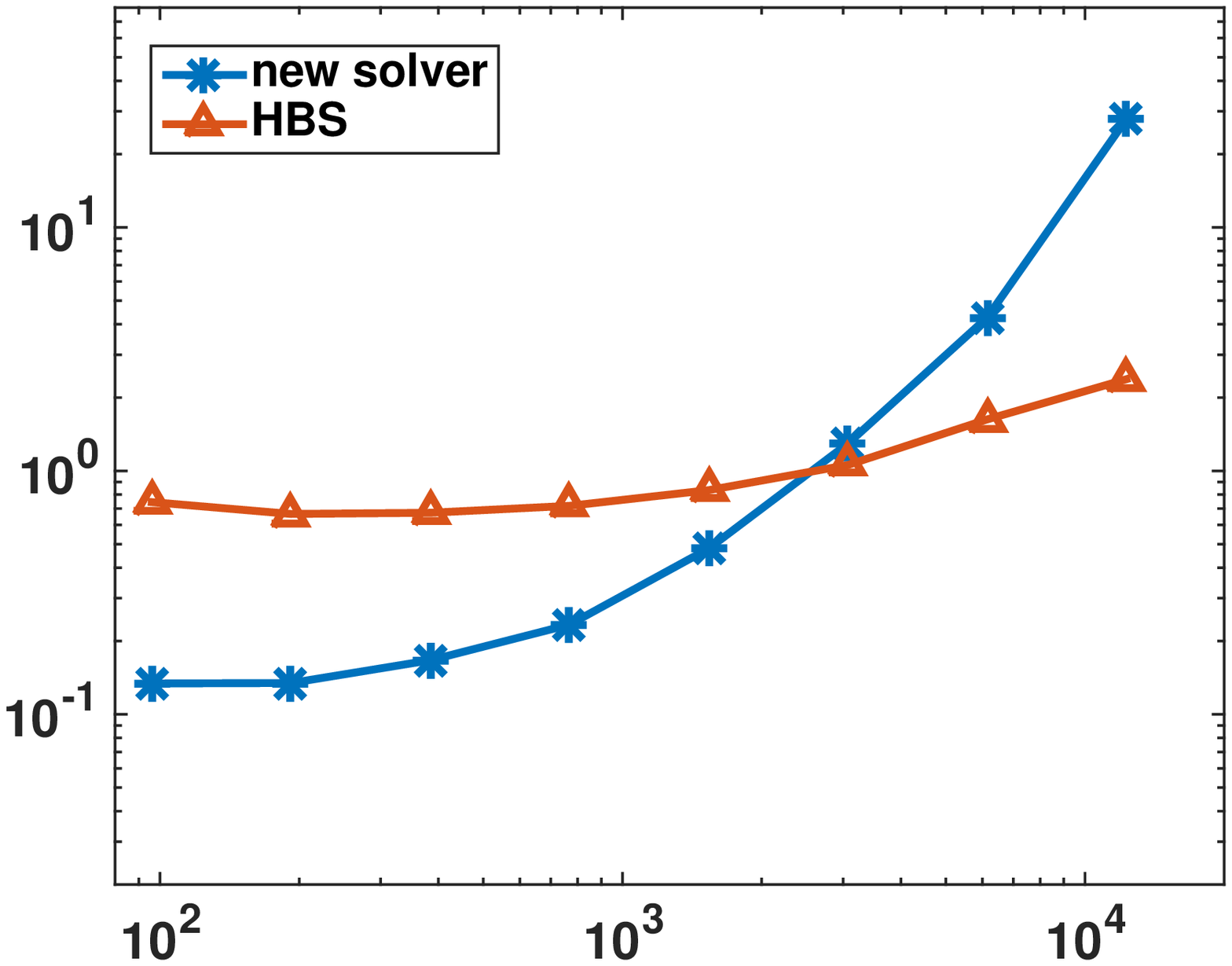}}
\put(200,06){\includegraphics[scale=0.4]{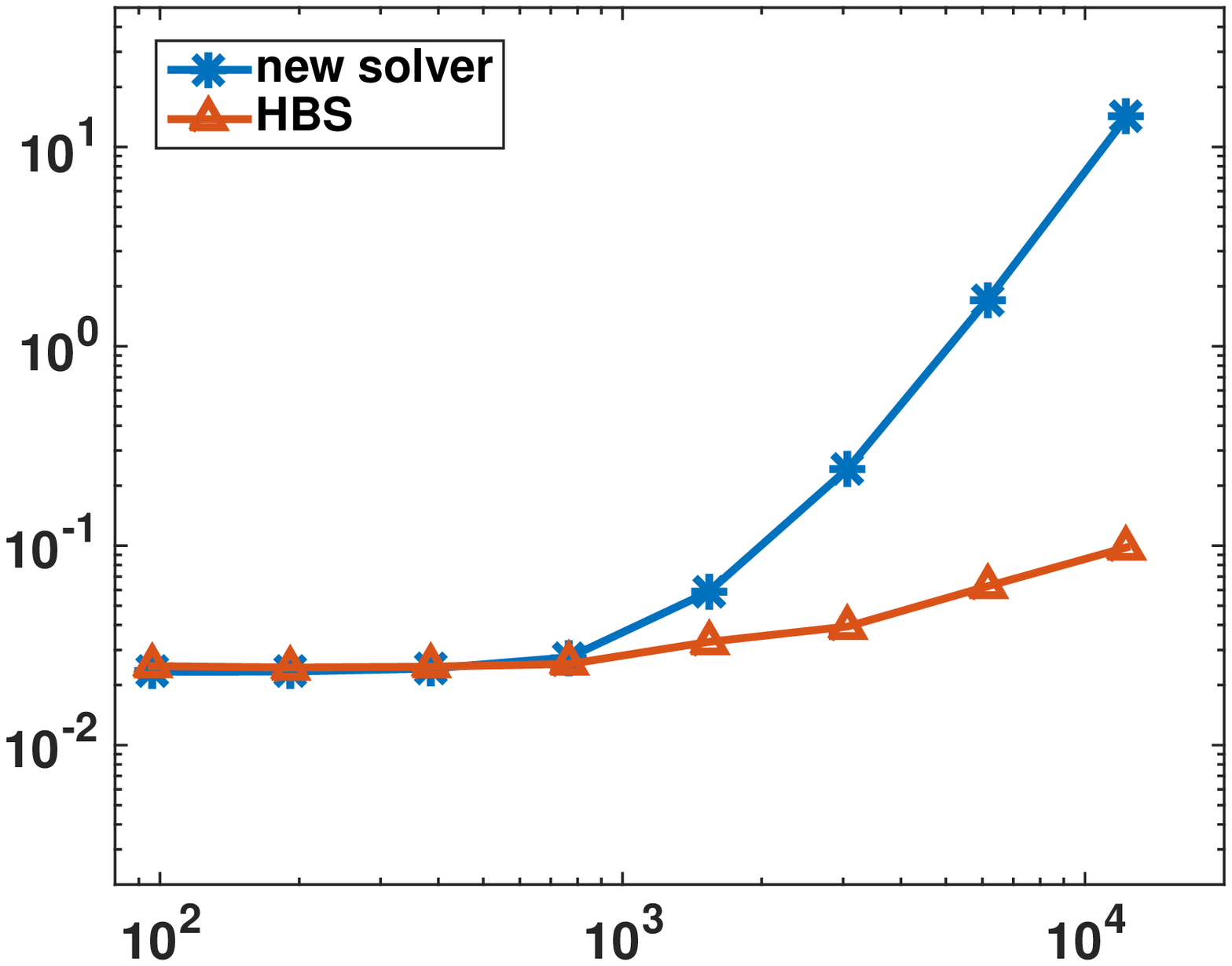}}
\put(50,175){ Precomputation}
\put(290,175){ Solve}
\put(-30,65){\rotatebox{90}{Time (sec)}}
\put(90,-10){$(a)$}
\put(300,-10){$(b)$}
\put(90,1){$N_p$}
\end{picture}

\end{center}
\caption{A log-log plot of the time in seconds versus $N_p$ for the 
(a) precomputation and (b) solve steps of both the new and HBS solver
for the star with refined panels geometry in Section \ref{sec:star_panel}. }
\label{fig:test2}

\end{figure}

This section considers the \textit{star with refined panels} geometry.  Such a geometry
occurs in the construction of an adaptive boundary integral equation discretization
technique. For this experiment $N_o$ is fixed, $N_o = 3200$, while $N_p$ increases.
Figure \ref{fig:test2} illustrates that even though the direct solver does not 
scale linearly with respect to $N_p$, there is a range of $N_p$ such that the 
new solver is faster than constructing a fast direct solver from scratch.  
In addition to the information listed in the beginning of this section, Table \ref{tab:test2} 
also reports the ratio of $N_p$ to $N_o$.  This ratio says that for $N_p$ less than half 
$N_o$ the new solver is at worst two times faster than building a new direct solver from 
scratch.  The addition of a fast direct solver applied $\mtx{A}_{pp}$ will
keep the speed up factor large for larger $N_p$.

\section{Summary}
\label{sec:summary}
This manuscript presented a fast direct solver for boundary value problems on 
locally perturbed geometries.  The solution technique is ideal for problems where
the local perturbation involves removing a small number of points.  Thus making the 
solver useful for optimal design problems where the perturbed geometry is placed in 
different portions of the original geometry and for improving the efficiency of  
adaptive refinement strategies.  For problems where the number of cut points $N_c$ is 
constant (corresponding to the optimal type problems) the method is three times faster 
than building a new direct solver from scratch.  For the adaptive refinement approach,
using the new solver is faster when the number of new points is less than fifty percent 
of the number of the original points on the geometry.

Future work will include the non-trivial extension of the solver presented in this 
paper to three dimensional boundary value problems.  This will involve careful 
management and processing of the geometry to make the best use of memory and 
limit communication.  The integration of the new solver to an adaptive integral 
equation discretization technique is also a future project.

\section{Acknowledgements} This research is supported by the Alfred P. Sloan foundation, the
NFS (DMS-1522631) and the Ken Kennedy Institute for Information Technology.

\bibliographystyle{abbrv} 
\bibliography{main_bib}

\end{document}